\documentclass[12pt]{article}

\usepackage{algorithm}
\usepackage{algpseudocode}
\usepackage{amsmath}
\usepackage{amssymb}
\usepackage{bm}
\usepackage{booktabs}
\usepackage{color}
\usepackage{cite}
\usepackage{epsfig}
\usepackage{epstopdf}
\usepackage{flafter}
\usepackage{float}
\usepackage{ifthen}
\usepackage{multirow}
\usepackage{amsfonts}
\usepackage{makecell} 
\usepackage{pifont}
\usepackage{subfigure}
\usepackage{times}
\usepackage{threeparttable}
\usepackage{url} 
\usepackage{subfigure,graphicx}
\usepackage{float}
\usepackage{color,multirow}
\usepackage{makecell}
\usepackage{setspace}
\usepackage{cite} 
\usepackage{url} 
\usepackage{ifthen} 
\usepackage[T1]{fontenc}
\usepackage{graphicx}
\usepackage{changepage}
\usepackage[colorlinks,linkcolor=black,anchorcolor=gray,citecolor=green,urlcolor=blue]{hyperref}
\usepackage{graphics}
\usepackage{booktabs}

\usepackage{amsfonts}
\usepackage{color}
\usepackage{amssymb}
\usepackage{amsmath}
\usepackage{indentfirst}

\def\sqr#1#2{{\vcenter{\vbox{\hrule height.#2pt
				\hbox{\vrule width.#2pt height#1pt \kern#1pt \vrule width.#2pt}
				\hrule height.#2pt}}}}
\def\signed #1{{\unskip\nobreak\hfil\penalty50
		\hskip2em\hbox{}\nobreak\hfil#1
		\parfillskip=0pt \finalhyphendemerits=0 \par}}
\def\endpf{\signed {$\sqr69$}}

\def\dbR{{\mathop{\rm l\negthinspace R}}}

\def\3n{\negthinspace \negthinspace \negthinspace }
\def\2n{\negthinspace \negthinspace }
\def\1n{\negthinspace }

\def\ds{\displaystyle}

\def\dbN{{\mathop{\rm l\negthinspace N}}}

\def\dbR{{\mathop{\rm l\negthinspace R}}}
\def\={\buildrel \triangle \over =}

\def\a{\alpha}

\def\d{\delta}
\def\e{\varepsilon}

\def\l{\lambda}

\def\n{\nabla}

\def\t{\times}

\def\th{\theta}
\def\o{\omega}

\def\i{\infty}
\def\ns{\noalign{\ss} }

\def\G{\Gamma}

\def\L{\Lambda}
\def\Si{\Sigma}

\def\O{\Omega}

\def\cC{{\cal C}}
\def\cD{{\cal D}}

\def\cK{{\cal K}}

\def\cO{{\cal O}}

\def\cQ{{\cal Q}}

\def\cS{{\cal S}}

\def\cl{{\cal l}}
\def\no{\noindent}

\def\ss{\smallskip}
\def\ms{\medskip}
\def\bs{\bigskip}
\def\q{\quad}
\def\qq{\qquad}

\def\max{\mathop{\rm max}}
\def\min{\mathop{\rm min}}
\def\exp{\mathop{\rm exp}}

\def\pa{\partial}

\def\cd{\cdot}

\def\div{\hbox{\rm div$\,$}}
\def\inf{\hbox{\rm inf$\,$}}

\def\cl{\overline}

\def\Re{{\mathop{\rm Re}\,}}

\def\|{\Big |}
\def\({\Big (}
\def\){\Big )}
\def\[{\Big[}
\def\]{\Big]}
\def\be{\begin{equation*}}
\def\bel{\begin{equation}\label}
\def\ee{\end{equation}}
\def\bt{\begin{theorem}}
\def\bcd{\begin{condition}}
\def\ecd{\end{condition}}
\def\et{\end{theorem}}
\def\bc{\begin{corollary}}
\def\ec{\end{corollary}}
\def\bde{\begin{definition}}
\def\ede{\end{definition}}
\def\bl{\begin{lemma}}
\def\el{\end{lemma}}
\def\bp{\begin{proposition}}
\def\ep{\end{proposition}}
\def\br{\begin{remark}}
\def\er{\end{remark}}
\def\ba{\begin{array}}
\def\ea{\end{array}}
\def\ed{\end{document}}
\def\ns{\noalign{\ms}}
\def\ds{\displaystyle}

\def\square#1{\vbox{\hrule\hbox{\vrule height#1%
			\kern#1\vrule}\hrule}}
\def\rectangle#1#2{\vbox{\hrule\hbox{\vrule height#1%
			\kern#2\vrule}\hrule}}

\font\tenbb=msbm10 \font\sevenbb=msbm7 \font\fivebb=msbm5

\newfam\bbfam
\scriptscriptfont\bbfam=\fivebb \textfont\bbfam=\tenbb
\scriptfont\bbfam=\sevenbb

\def\oO{{\overline \O}}

\newtheorem{lemma}{Lemma}[section]
\newtheorem{remark}{Remark}[section]

\newtheorem{theorem}{Theorem}[section]
\newtheorem{corollary}{Corollary}[section]

\newtheorem{definition}{Definition}[section]
\newtheorem{proposition}{Proposition}[section]
\newtheorem{condition}{Condition}[section]

\makeatletter

\@addtoreset{equation}{section}
\makeatother
\begin{document}

\title{\bf Observability estimate for the  wave equation with variable coefficients \thanks{This work is partially supported by the NSF of China
		under grant 11971333, 11931011, and by the Science Development Project of Sichuan University under grant 2020SCUNL201. The authors gratefully acknowledge Professor Xu Zhang for his help and encouragement.}}
\author{Xiaoyu Fu\thanks{School of Mathematics,  Sichuan University,
		Chengdu 610064, China. ({\tt xiaoyufu@scu.edu.cn}).} \and  Zhonghua Liao\thanks{School of Mathematics, Sichuan University, Chengdu 610064, China. ({\tt zhonghualiao@yeah.net}).}}

\date{}

\maketitle

\begin{abstract}This paper is devoted to a study of  observability estimate for the wave equation with variable coefficients $(h^{jk}(x))_{n\t n}$ ($n\in\dbN)$. We consider both the observation point lies outside the domain and the observation point lies inside the domain. Based on  a Carleman estimate for the ultra-hyperbolic operator and a delicate treatment of observation region, we obtain two observability estimates with explicit observability constants.   The key improvements are: (1) we improve the requirement of waiting time $T$; (2) we improve the size of the observation region (see Fingure \ref{Fig1} and Fingure \ref{Fig0} for the case of $(h^{jk}(x))_{n\t n}=I_n)$.

\end{abstract}

\no{\bf 2010 Mathematics Subject Classification}. 93B05, 93B07, 93B27

\bs

\no{\bf Key Words}. Observability estimate,  Hyperbolic equation, Carleman estimate.
\pagestyle{plain}\thispagestyle{plain}
\section{Introduction}
Given $ T>0 $ and a bounded domain $ \O $ of $ \dbR^n (n\in\dbN)$ with $ C^2 $ boundary $ \Gamma $, put $ Q=(0,T)\t\O, ~ \Sigma=(0,T)\t\Gamma. $

 Let  $ h^{jk}(\cd)\in C^2(\cl{\O}) $  be fixed satisfying
\begin{equation}\label{a1}
h^{jk}(x)=h^{kj}(x),\q\q \forall x\in \cl{\O},~~\mbox{ }  j,~k=1,\cd\cd\cd, n,
\end{equation}
and for some constant $ h_0>0 $,
\begin{equation}\label{a2}
\sum_{j,k=1}^n h^{jk}(x)\xi^j\bar{\xi}^k\ge h_0|\xi|^2,\q\q \forall (x,\xi)\in \cl{\O}\t\dbR^n,\q\xi=(\xi^1,\cd\cd\cd,\xi^n).
\end{equation}
 Let us consider the following wave equation with variable coefficients:
\bel{aaa1}
\left\{
\ba{ll}
\ds w_{tt}-\sum_{j,k=1}^n (h^{jk}w_{x_j})_{x_k}=qw+\sum_{k=1}^n q_1^k w_{x_k}+q_2w_t,\q &\mbox{in } Q,\\
\ns\ds w=0,&\mbox{on } \Si,\\
\ns\ds w(0)=w_0,\q w_t(0)=w_1,& \mbox{in } \O,
\ea
\right.
\ee
where $(w_0,w_1)\in L^2(\O)\t H^{-1}(\O)$, $q\in L^\infty (Q) $ and $ q_1^k\in W^{1,\i} (Q) $ ($ k=1,\cd\cd\cd,n $), $ q_2\in W^{1,\i}(Q) $. By the method of (\cite[Ch. I, Th 4.2]{Lions}) it follows that the existence and the uniqueness of the solution $w$ of (\ref{aaa1})  lies in the class $\ds w\in C([0,T]; L^2(\O))\cap C^1([0,T]; H^{-1}(\O))$.

The main purpose of this paper is to study the  internal observability problem of (\ref{aaa1}), by which we mean the following: given $ T>0 $ and $ \cK $ be a sub-domain of $ Q $, find (if possible) a constant $ \cC=\cC(q, \{q_1^k\}_{k=1}^n, q_2)>0 $ such that the corresponding solution $w$ of (\ref{aaa1}) satisfies
\begin{equation}\label{abc1}
|w_0|_{L^2(\O)}^2+|w_1|_{H^{-1}(\O)}^2\le \cC(q, \{q_1^k\}_{k=1}^n, q_2)\int_{\cK}|w|^2dxdt, \q \forall (w_0,w_1)\in L^2(\O)\t H^{-1}(\O).
\end{equation}

The inequality (\ref{abc1}) is called an observability estimate for (\ref{aaa1}). This inequality means that the initial energy of a solution in the time $ t=0 $ can be bounded by its partial energy in the local observation region $\cK $. Such kind of inequalities are closely related to control and state observation problems of wave equations with constant or variable coefficients. For example, they can be applied to a study of the controllability (e.g. \cite{BLR, JMC, XYZ, Lions, XZ2, ZZ, EZ1, DR}), the stabilization of some locally damped semilinear wave equation with variable coefficients (e.g. \cite{LT}) and also inverse problems (e.g. \cite{MVK, MML}). In this respect, there exist numerous works devoted to observability estimates of wave equations with constant or variable coefficients, we refer to \cite{DZZ, LITR, XQX, Shao, XZ1, XZ4, JLQL} and rich references therein.

Note that the constant $ \cC $ in (\ref{abc1}) depends on the lower-order term coefficients $ q(\cd) ,\ q_1^k(\cd) (k=1,2, \cdots, n), q_2(\cd)$, the observation domain $ \cK $ and the waiting time $ T $ in (\ref{aaa1}). In this paper, the explicit estimate of $C$, the critical value of  waiting time $T$, the size of the observation region $\cK$ are parts of the problems we concerned.

For the case  $ (h^{jk}(\cd))_{n\t n}=I_n $, $ ~q(\cd)\in L^\i(Q) $, $ q_1^k(\cd)=0\  (k=1,2,\cdots, n)$, $ q_2(\cd)=0 $, it was proved  in Zhang \cite{XZ1} that explicit observability estimate (\ref{abc1}) holds  for
\bel{1201-a}
T> 2\max_{x\in\O} |x-x_0| ,\q\cK=\cK_1\= (0,T)\t\o,\q C(q)=C\exp(\exp(\exp(C|q(\cd)|_{L^\i(Q)}))),
\ee
where $ x_0\in \dbR^n\backslash \oO $, $ \o $ is some given neighborhood of $ \G_0 $, and $ \G_0 $ is part of  the boundary of $ \G $ satisfying certain conditions, which will be specified later. For the case $ (h^{jk}(\cd))_{n\t n}=I_n,~q(\cd)\in C^\i(Q) $, $ \{q_1^k(\cd)\}_{k=1}^n\in C^\i(Q) $, $ q_2(\cd)\in C^\i (Q)$, it was proved in Jena \cite{VJ} that observability estimate (\ref{abc1}) holds for
$$
x_0\in \dbR^n\backslash \oO ,\q T> 2\max_{x\in\O} |x-x_0| ,\q\cK=\cK_2\= \cK_1\cap \Big\{(t,x)\in \dbR^{1+n}\big|~|x-x_0|^2>t^2\Big\},
$$
and
$$
x_0\in  \oO,\q  T> 2\max_{x\in\O} |x-x_0| ,\q\cK=\cK_3,
$$
where $\cK_3 $ is some neighborhood of  $\cK_2$. We also refer to the related references \cite{Shao} for boundary observability estimate of linear wave equation (i.e. $ (h^{jk}(\cd))_{n\t n}=I_n $) on time-dependent domains with a smaller boundary observation region. In the case of wave equation with variable coefficients, $ ~q(\cd)\in L^\i(0, T; L^n(\O)) $, $ q_1^k(\cd)=0\  (k=1,2,\cdots, n)$, $ q_2(\cd)=0 $, it was proved  in  \cite{XYZ} that explicit observability estimate (\ref{abc1}) holds  for  larger waiting time   and observation region $\cK=\cK_1$ satisfying (\ref{1201-a}).  In this paper, we consider the explicit observability inequality of (\ref{aaa1}) with variable coefficients. Based on a Carleman estimate for the ultra-hyperbolic operator and a delicate treatment on the region of observation, we improve not only the requirement of waiting time $ T $ but also the size of the observation region (see Theorem \ref{theorem1} and Theorem \ref{theorem2} in detail). Particularly, for the case $ (h^{jk}(\cd))_{n\t n}=I_n $, the waiting time  $ T $ can be improved to $\ds  T>2\max_{x\in \oO\backslash \o}|x-x_0| $, the changes of the observation regions can be seen in Remark \ref{rem}, i.e. Figure \ref{Fig1} and Figure \ref{Fig0}.

The rest of this paper is organized as follows. In Section 2, we give the statement of our main results. In Section 3, we collect some preliminaries we needed. Finally, we give the proofs of our main results in Section 4 and Section 5, respectively.

\section{Statement of the main results}

To begin with, we introduce the following condition:

\bcd\label{cd1}
There exists a function  $d(\cd)\in C^2(\oO)$ satisfying the
following:

\ss

{\rm(i)} For some constant $\mu_0>0$, it holds
\bel{a8}
\ba{ll}\ds
\sum_{j,k=1}^n\Big\{\sum_{j', k'=1}^n\[2h^{jk'}
(h^{j'k}d_{x_{j'}})_{x_{k'}}-h^{jk}_{x_{k'}}h^{j'k'}d_{x_{j'}}
\]\Big\}\xi^j\xi^k\geq \mu_0\sum_{j,k=1}^n h^{jk}\xi^j\xi^k,\\
\ns
\ds
\qq\qq\qq\qq\qq\qq\qq\forall\; (x,\xi^1,\cdots,\xi^n)\in\oO\times
\dbR^n.
\ea
\ee

{\rm(ii)} There is no critical point of function $d(\cd)$ in $\oO$,
i.e.,
\bel{1201-a6}
\min_{x\in\oO} |\n d(x)|> 0.
\ee
%
\ecd

Note that for the case $ (h^{jk})_{n\t n}=I_n $, and any given $ x_0\in \dbR^n\backslash \oO $, by choosing $ d(x)=|x-x_0|^2 $, we have (\ref{a8}) with $ \mu_0=4 $. In this case,
\begin{equation}\label{abc5}
\G_0=\Bigg\{x\in \G\| ~(x-x_0)\cd \nu(x)>0\Bigg\}.
\end{equation}
We refer to \cite{XYZ, YL} for  more examples and explanations on  Condition \ref{cd1}.

For the function $d(\cd)$ satisfying Condition \ref{cd1}, we
introduce the following set:

\begin{equation}\label{abc4}
\G_0\=\left\{x\in \G\| ~\sum_{j,k=1}^n h^{jk}d_{x_j}\nu_{k}>0\right\}.
\end{equation}

It is easy to check that, if $ d_0(\cd)\in C^2(\cl{\O}) $ satisfies (\ref{a8}), then for any given constants $ a\ge 1 $  and $ b\in \dbR $, the function
\begin{equation}\label{a14}
\hat{d}=\hat{d}(x)\=ad_0(x)+b
\end{equation}
still satisfies Condition \ref{cd1} with $ \mu_0 $ replaced by $ a\mu_0 $, meanwhile, the scaling and translating $ d_0(x) $ do not change the set $ \G_0 $. Hence, by scaling and translating $ d(x) $, if necessary, we may assume without loss of generality that
\begin{equation}\label{a15}
\left\{
\ba{ll}
\ns\ds (\ref{a8}) \mbox{ holds with } \mu_0\ge 4,\\
\ns\ds \frac{1}{4}\sum_{j,k=1}^nh^{jk}(x){d}_{x_j}{d}_{x_k}\ge {d}(x)>0, \q\q \forall x\in \cl{\O}.
\ea
\right.
\end{equation}

For any set $ M\subset\dbR^n $ and $ \d>0 $, we define $ \cO_\d(M)=\big\{x\in\dbR^n|~|x-x'|<\d \mbox{
	for  } x'\in M\big\} $. Let  $\o, \o_0$ be proper open
non-empty subsets of $\O$ satisfying $\o_0\subset\subset \o $. We assume that  there exist constants  $ 0< \d_0<\d $ such that
\begin{equation}\label{abc6}
\o\=\cO_\d (\G_0)\bigcap \O, \q \o_0 \=\cO_{\d_0}(\G_0)\bigcap \O, \q Q_{\o\backslash\o_0}\=(0,T)\t \{\o\backslash \o_0\}.
\end{equation}
In what follows, we set
\begin{equation}\label{abc7}
R_0\=\min_{x\in \oO}\sqrt{d(x)},\q R_1\= \max_{x\in \oO\backslash \o} \sqrt{d(x)}, \q T_*\= 2\inf\Big\{R_1\big| ~d(\cd) \mbox{ satisfies } (\ref{a15})\Big\}.
\end{equation}
For the function $ d(\cd) $ satisfying (\ref{abc7}), for any constant $ \d_1 \in (0, 1/2) $, we define
\begin{equation}\label{abc8}\left\{\ba{ll}\ds
\cD \=\Big\{(t,x)\in Q\big|~ d(x)-(t-T/2)^2>0\Big\},\\
\ns\ds
\cK\= \Big\{\(\frac{T}{2}-\d_1 T, \frac{T}{2}+\d_1 T\)\t \o_0\Big\}\cup \left\{Q_{\o\backslash\o_0} \cap \cD\right\}.
\ea\right.\end{equation}

Throughout of this paper, we shall denote by $ C=C(\O, n, (h^{jk})_{n\t n},T) $ a generic positive constant, which may change from line to line (unless otherwise stated). We have the following observability result.
\begin{theorem}\label{theorem1}
Let $ h^{jk}(\cd)\in C^2(\oO) $ satisfy (\ref{a1}) and (\ref{a2}), $ q(\cd)\in L^\i (Q) $ and $ q_1^k(\cd)\in W^{1,\i} (Q) $ ($ k=1,\cd\cd\cd,n $), $ q_2(\cd)\in W^{1,\i}(Q) $. Let the function $d(\cd)$ satisfies Condition \ref{cd1} and (\ref{a15}). Then for any $ T>T_* $, the weak solution $ w(\cd)\in C([0,T];L^2(\O))\cap C^1([0,T];H^{-1}(\O)) $ of equation (\ref{aaa1}) satisfies that
\begin{equation}\label{abc10}
|w_0|_{L^2(\O)}^2+|w_1|^2_{H^{-1}(\O)}\le \cC(r) \int_\cK |w|^2 dxdt, \q \forall (w_0, w_1)\in L^2(\O)\t H^{-1}(\O),
\end{equation}
where
\begin{equation}\label{aabb1}\left\{\ba{ll}\ds
\cC(r)=C\exp (\exp (Cr)),\\
\ns\ds r\= \max\Big\{|q(\cd)|_{L^\i(Q)}, \{|q_1^k(\cd)|_{W^{1,\i}(Q)}\}_{k=0}^n,|q_2(\cd)|_{W^{1,\i} (Q)} \Big\}.
\ea\right.\end{equation}

\end{theorem}

\begin{figure}[!h]
\footnotesize
\begin{center}
\begin{tabular}{ccc}
\includegraphics[width=0.2\textwidth]{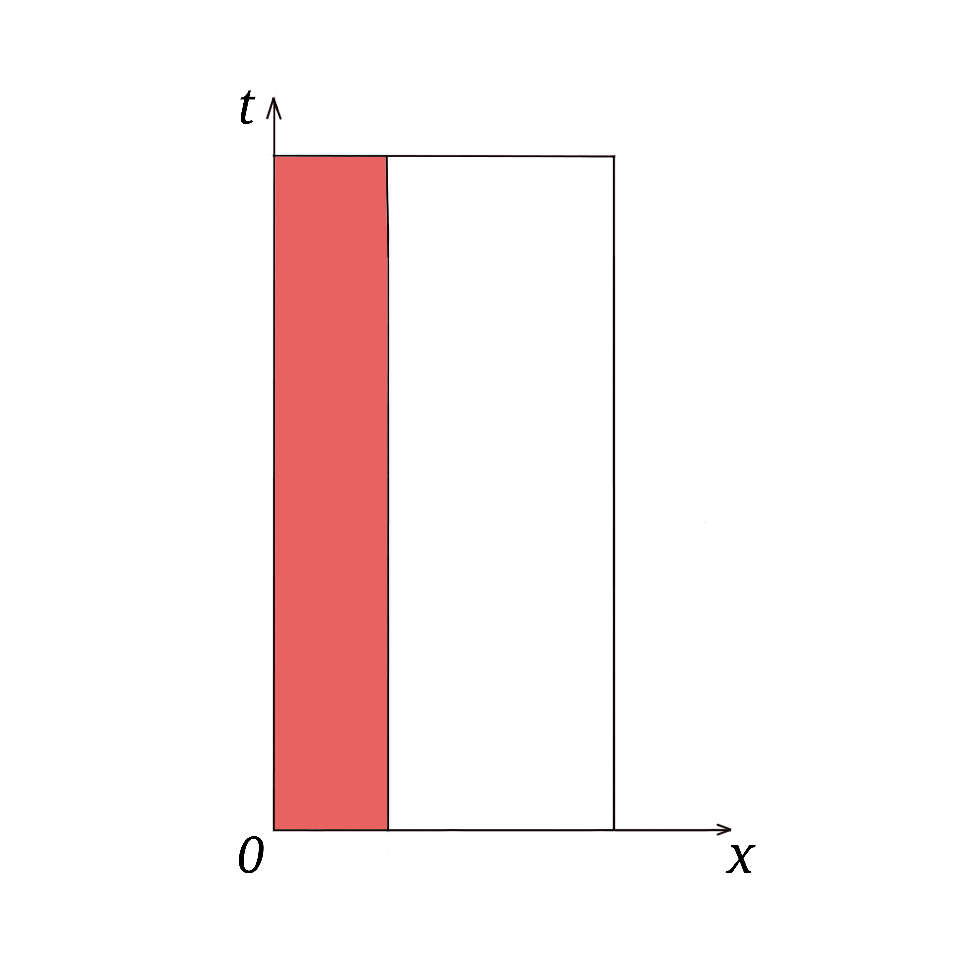}&
\hspace{-0.46cm}
\includegraphics[width=0.2\textwidth]{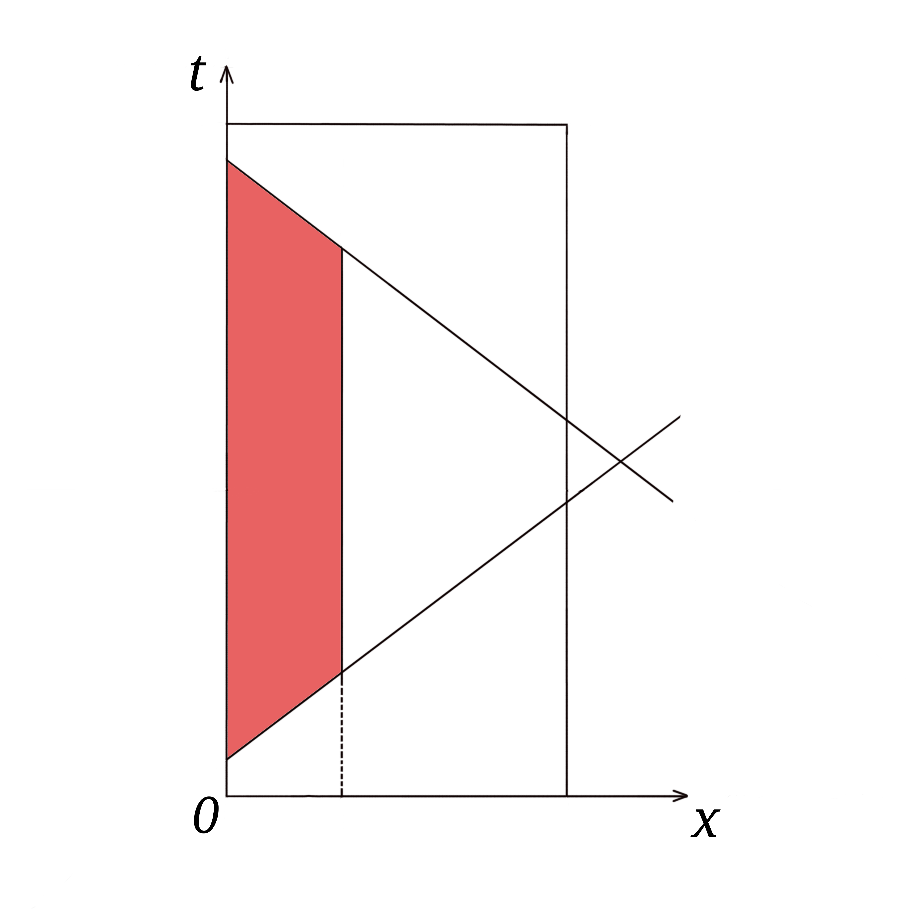}&
\hspace{-0.46cm}
\includegraphics[width=0.2\textwidth]{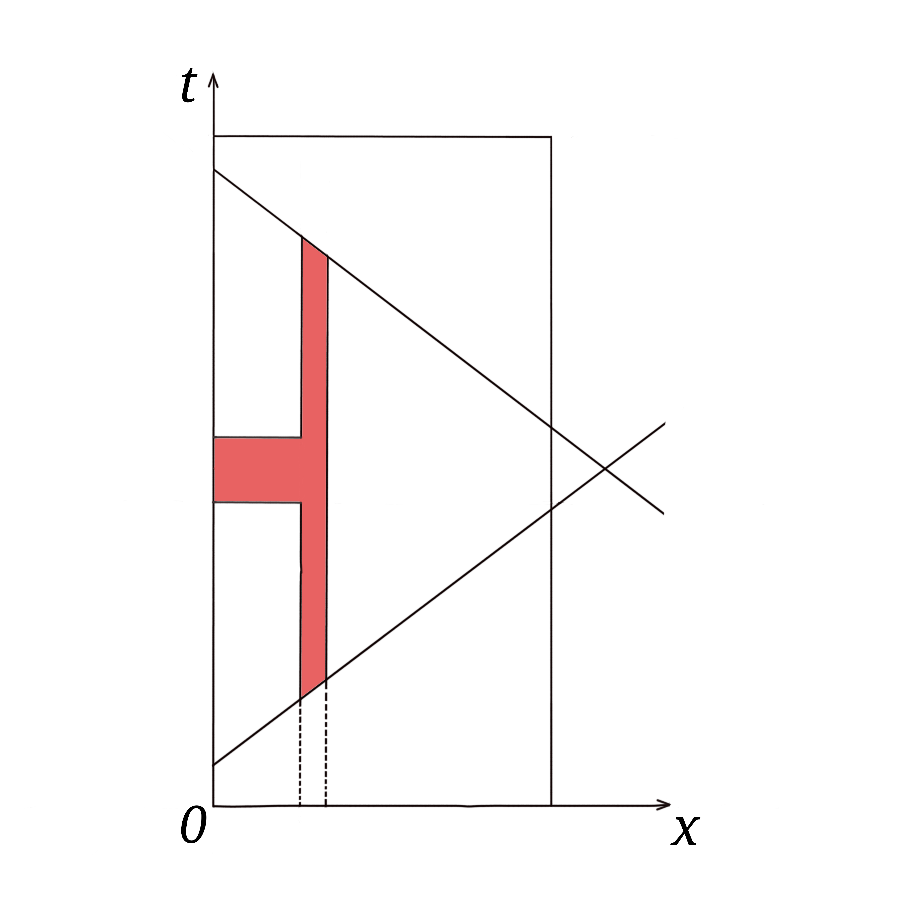}\\
(a)&(b)&(c)
\hspace{-0.46cm}
\end{tabular}

\end{center}
\caption{For the wave equation in one dimensional case, the red regions in (a) and  (b) are the observation regions developed in \cite[Theorem 2.1]{XZ1} and \cite[Theorem 1.2]{VJ}, respectively. The red region in (c) is our result stated in Theorem \ref{theorem1}. }\label{Fig1}
\end{figure}

\begin{figure}[!h]
\centering
\includegraphics[width=3.5in]{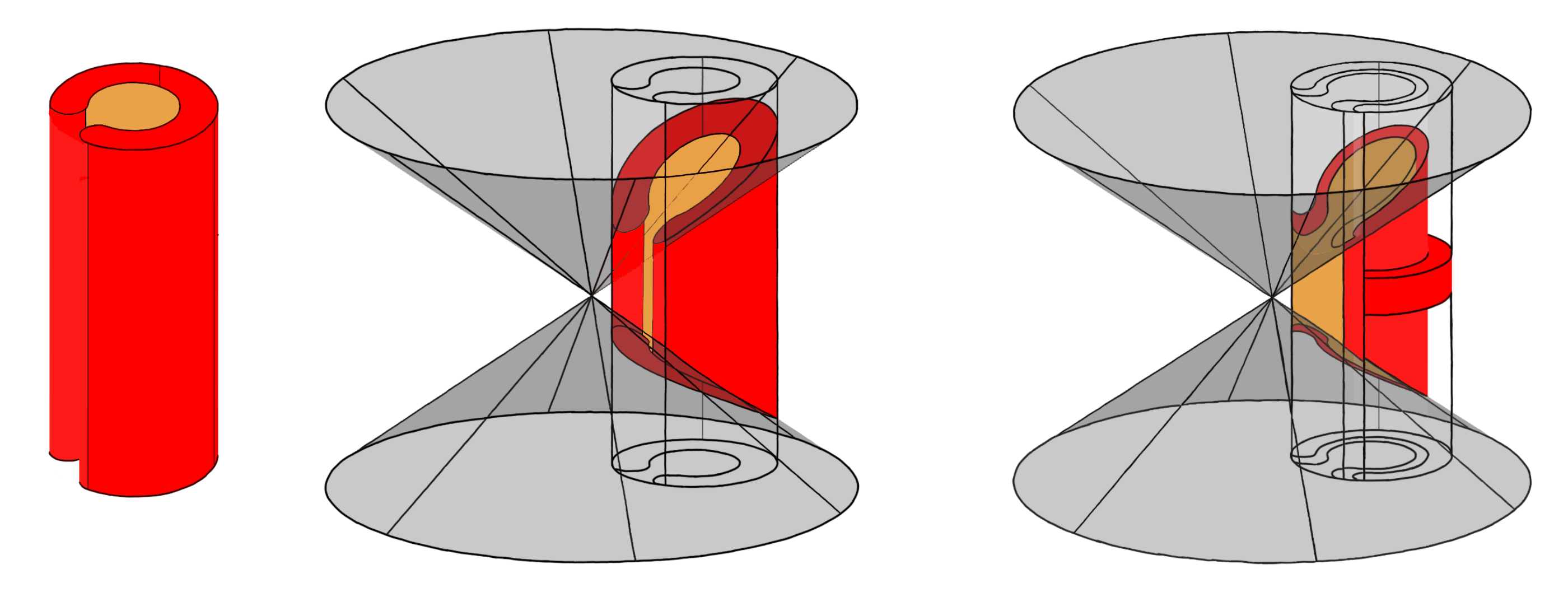}
\caption{The changes of the observation regions for the wave equations in multidimensional case.
}\label{Fig0}
\end{figure}

\begin{remark}\label{rem}
In the case of $\ds (h^{jk})_{n\t n}=I_n$, we refer to Figure \ref{Fig1} for the  changes of observation regions in one dimensional case, and Figure \ref{Fig0} for the changes of observation regions in multidimensional case.
\end{remark}

Note that Theorem \ref{theorem1} holds provided the function $d(\cd)$ no critical point in $\oO$,
i.e., $\ds \min_{x\in\oO} |\n d(x)|> 0 $.  In the case that there is a critical point of function $d(\cd)$ in $\oO$, we  can  also establish suitable observability estimate.  To this aim, we introduce the following Condition.

\begin{condition}\label{condition2}
There exists a function $ b(\cd)\in C^2(\oO) $ satisfying the following:

\ss

{\rm(i)} For some constant $\mu_0>0$, (\ref{a8}) holds for $b(\cd)$.

{\rm(ii)} There is only one critical point $ x_0\in \oO$, i.e.
\begin{equation}\label{abc11}
\min_{x\in \oO\backslash \{x_0\}}|\n b(x)|>0, \q |\n b(x_0)|=0.
\end{equation}

{{\rm(iii)} There exists constant $s>0$ such that
	\begin{equation}\label{abc12}
	\min_{x\in \oO}b(x)=b(x_0)=0, \q\lim_{x\to x_0}\frac{\sum_{j,k=1}^n h^{jk}(x)b_{x_j}(x)b_{x_k}(x)}{b(x)}=s.
	\end{equation}}
\end{condition}
\par Let $\zeta$ ($|\zeta|\neq  0$) be sufficient small vector in $\dbR^n$ satisfying $x_0-\zeta\in \oO$ and there exist constant $0<\d_2<\d_0$, such that $\o_0\supset \cO_{\d_2} \(\G_{0,\zeta}\)\cap \O$ where
\bel{gab1}
\G_{0,\zeta}\= \Big\{ x\in \G\|~\sum_{j,k=1}^n h^{jk}d_{x_j}(x+\zeta)\nu(x)>0\Big\}.
\ee
We have that if function $d(\cd), d(\cd+\zeta)\in C^2(\oO)$ satisfy  Condition \ref{condition2},
we can also establish suitable observability estimate. Similar to the analysis of  (\ref{a15}),  without loss generality, we assume the following hold:
\begin{equation}\label{a15-a}
\left\{
\ba{ll}
\ns\ds (\ref{a8}) \mbox{ holds with } \mu_0> 4,\\
\ns\ds \frac{1}{4}\sum_{j,k=1}^nh^{jk}(x){d}_{x_j}(x){d}_{x_k}(x) > d(x)>0, \q\q \forall x\in \O\backslash \{x_0\},\\
\ns\ds \frac{1}{4}\sum_{j,k=1}^n h^{jk}(x)d_{x_j}(x+\zeta)d_{x_k}(x+\zeta)>d(x+\zeta)>0,\qq \forall x\in \O\backslash \{x_0-\zeta\}.
\ea
\right.
\end{equation}
In what follow, we set
\bel{lzh1}
R_1\= \max_{x\in \oO\backslash \o}\Big\{ \sqrt{d(x)}, \sqrt{d(x+\zeta)}\Big\},\q T^*\=2\inf\Big\{ R_1\|~d(\cd) \mbox{ satisfies } (\ref{a15-a})\Big\},
\ee
and
\bel{ga1}
\ba{ll}\ds
\cD_\zeta\=\Big\{(t,x)\in Q\|~ d(x+\zeta)-(t-T/2)^2>0\Big\},\\
\ns\ds \cK_\zeta\=\Big\{\(\frac{T}{2}-\e T, \frac{T}{2}+\e T\)\t \o_0\Big\}\cup \{Q_{\o\backslash\o_0}\cap \cD_\zeta\}.
\ea
\ee
We have the following observability result.
\begin{theorem}\label{theorem2}
Let $ h^{jk}(\cd)\in C^2(\oO) $ satisfy (\ref{a1}) and (\ref{a2}), $ q(\cd)\in L^\i (Q) $ and $ q_1^k(\cd)\in W^{1,\i} (Q) $ ($ k=1,\cd\cd\cd,n $), $ q_2(\cd)\in W^{1,\i}(Q) $. Let $d(\cd), d(\cd+\zeta)$ hold Condition \ref{condition2} and (\ref{a15-a}),  $\cK$  and $\cK_\zeta $ be given by (\ref{abc8}) and (\ref{ga1}). For any open domain $W$ satisfies
$\ds
\overline{\cK}\cup \cl{\cK_\zeta}\subset W\subset (-T,T)\t \O,
~ T>T^* $, the weak solution $ w(\cd)\in C([0,T];L^2(\O))\cap C^1([0,T];H^{-1}(\O)) $ of equation (\ref{aaa1}) satisfies
\begin{equation}\label{abc14}
|w_0|_{L^2(\O)}^2+|w_1|^2_{H^{-1}(\O)}\le \cC(r) \int_W |w|^2 dxdt, \q \forall (w_0, w_1)\in L^2(\O)\t H^{-1}(\O),
\end{equation}
where $ \cC(r) $ is the same as (\ref{aabb1}).
\end{theorem}
\begin{remark}
For the case $(h^{jk})_{n\t n}=I_n$, by choosing $d(x)=|x-x_0|^2$ ($x_0$ inside $\oO$), for any open domain $(0,T)\t \O\supset  W \supset \cl{\cK} $ and $T>T^*$, there exist sufficient small $\zeta$ with $x_0-\zeta\in \oO$ such that  $W\supset \cl{\cK}\cup \cl{\cK_\zeta}$ and $T>\tilde T^*$ where
\bel{ha1}
\tilde T^*\= 2\max_{x\in \oO\backslash \o} |x-x_0|.
\ee
Notice that $|x-x_0|^2$ and $|x+\zeta-x_0|^2$ satisfy Condition \ref{condition2} and (\ref{a8}), so that (\ref{abc14}) hold with $W$.
\end{remark}

\section{Some preliminaries}
In this section, we collect some preliminaries we needed. First, based on a weighted identity established in \cite{XQX}, we have the following pointwise inequality for the ultra-hyperbolic operator $\ds  \pa_t^2+\pa_s^2-\sum_{j,k=1}^n\pa_{x_k}(h^{jk}x_{x_k}) $:
\begin{corollary}\label{Co4.1}
Let $ u\in C^2(\dbR^{2+n};\dbR) $, $ \ell\in C^3(\dbR^{2+n}; \dbR) $ and $ \Psi \in C^1( \dbR^n; \dbR) $. Set $ \th=e^\ell $ and $ v=\th u $. Then,
\begin{equation}\label{c1}
\ba{ll}
\ds \th^2 \|u_{tt}+u_{ss}-\sum_{j,k=1}^n(h^{jk}u_{x_j})_{x_k}\|^2+2\div V+2M_t+2N_s\\
\ds \ge 2\[\ell_{tt}-\ell_{ss}+\sum_{j,k=1}^n(h^{jk}\ell_{x_j})_{x_k}+\Psi\]v_t^2-8\sum_{j,k=1}^nh^{jk}\ell_{tx_j}v_{x_k}v_t+8\ell_{st}v_sv_t\\
\ds\q-8\sum_{j,k=1}^nh^{jk}\ell_{sx_j}v_{x_k}v_s +2\[\ell_{ss}-\ell_{tt}+\sum_{j,k=1}^n(h^{jk}\ell_{x_j})_{x_k}+\Psi\]v_s^2\\
\ds\q +2\sum_{j,k=1}^n c^{jk}v_{x_j}v_{x_k}-2\sum_{j,k=1}^nh^{jk}\Psi_{x_j}vv_{x_k}+Bv^2,
\ea
\end{equation}
where
\begin{equation}\label{c2}
\left\{\ba{ll}
\ds A=\sum_{j,k=1}^n (h^{jk}\ell_{x_j}\ell_{x_k}-h_{x_j}^{jk}\ell_{x_k}-h^{jk}\ell_{x_jx_k})-\ell_t^2-\ell_s^2+\ell_{tt}+\ell_{ss}-\Psi,\\
\ns\ds c^{jk}=\sum_{j',k'=1}^n\[2h^{jk'}(h^{j'k}\ell_{x_{j'}})_{x_{k'}}-(h^{jk}h^{j'k
	'}\ell_{x_{j'}})_{x_{k'}}\]+h^{jk}(\ell_{tt}+\ell_{ss}-\Psi ),\\
\ns\ds B =2\[A\Psi -(A\ell_t)_t-(A\ell_s)_s+\sum_{j,k=1}^n(Ah^{jk}\ell_{x_j})_{x_k}\],
\ea\right.
\end{equation}
and
\begin{equation}\label{c3}
\left\{\ba{ll}
\ns\ds V=[V^1,\cd\cd\cd, V^k,\cd\cd\cd , V^n],\\
\ns\ds V^k=2\sum_{j,j',k'=1}^nh^{jk}h^{j'k'}\ell_{x_{j'}}v_{x_j}v_{x_{k'}}+\sum_{j=1}^nh^{jk}A\ell_{x_j}v^2-\Psi v\sum_{j=1}^n h^{jk}v_{x_j}\\
\ns\ds \q\q \q-\sum_{j,j',k'=1}^nh^{jk}h^{j'k'}\ell_{x_j}v_{x_{j'}}v_{x_{k'}}-2(\ell_tv_t+\ell_sv_s)\sum_{j=1}^n h^{jk}v_{x_j}\\
\ns\ds \qq\q+\sum_{j=1}^n h^{jk}\ell_{x_j} (v_t^2+v_s^2),\\
\ns\ds M_t=\ell_t \(v_t^2-v_s^2+\sum_{j,k=1}^nh^{jk}v_{x_j}v_{x_k}\)-2\sum_{j,k=1}^nh^{jk}\ell_{x_j}v_{x_k}v_t+2\ell_s v_s v_t\\
\ns\ds \qq\q+\Psi vv_t-A\ell_t v^2,\\
\ns\ds N_s=\ell_s \(v_s^2-v_t^2+\sum_{j,k=1}^nh^{jk}v_{x_j}v_{x_k}\)-2\sum_{j,k=1}^nh^{jk}\ell_{x_j}v_{x_k}v_s+2 \ell_t v_s v_t\\
\ns\ds \qq\q+\Psi vv_s-A\ell_s v^2.
\ea\right.
\end{equation}
\end{corollary}
\par {\it Proof.}  In \cite[ Theorem 1.1]{XQX},  we choose
\begin{equation}\label{abc18}
(a^{jk})_{(n+2)\t (n+2)}=\left(\ba{ccc}
-1&0&0\\
0&-1&0\\
0&0&(h^{jk})_{n\t n}
\ea\right).
\end{equation}
By elementary computation, one can obtain Corollary \ref{Co4.1} immediately. For the reader's convenience,  we also give a detailed proof of Corollary \ref{Co4.1} in Appendix A.~\endpf

Next, we recall the following known result.
\begin{lemma}\label{lll}
{\rm (\cite{Evans} )} Denote $ L\= \sum_{j,k=1}^n \pa_{x_j}(h^{jk}\pa_{x_k}) $. There exist a constant $ \L $, such that for all $ \Re \l \ge \L $, equation $ (-L+\l)u=f $ have unique solution $ u\in H_0^1(\O) $ for all $ f\in H^{-1}(\O) $ and there exist  positive constants $ C_1 $ and $C_2$ such that
\begin{equation}\label{bca3}
C_1|u|_{H_0^1(\O)}\le |f|_{H^{-1}(\O)}\le C_2 |u|_{H_0^1(\O)}.
\end{equation}
\end{lemma}
Put
\begin{equation}\label{abc19}
E(t)\= \frac{1}{2}\[|w_t(t,\cd)|^2_{H^{-1}(\O)}+|w(t,\cd)|^2_{L^2(\O)}\].
\end{equation}
We have the following energy estimate.
\begin{lemma}\label{lemma31}
Let $ T>0 $, $q(\cd)\in L^\infty (Q) $ and $ q_1^k(\cd)\in W^{1,\i} (Q) $ ( $ k=1,\cd\cd\cd,n $) , $ q_2(\cd)\in W^{1,\i}(Q) $, $ w_0(\cd)\in L^2(\O) $ and $ w_1(\cd)\in H^{-1}(\O) $. Then the weak solution $ w(\cd)\in C([0,T]; L^2(\O))\cap C^1([0,T]; H^{-1}(\O)) $ of (\ref{aaa1}) satisfies
\begin{equation}\label{abc20}
E(t)\le CE(s)e^{Cr},\qq \forall t,s\in [0,T],
\end{equation}
where $r$ is given by (\ref{aabb1}).
\end{lemma}

{\it Proof of Lemma \ref{lemma31}. } By Lemma \ref{lll}, we know that there exists a constant $ \l_0 $, such that  operator $ -L+\l_0 $ is a reversible operator. Put $f\= (-L+\l_0)^{-1} w_t$, we have
\bel{ha2}
\ba{ll}
\ds w_t(-L+\l_0)^{-1}\( w_{tt}+(-L+\l_0)w\)\\
\ns\ds =(-L+\l_0) f f_t+w_t w\\
\ns\ds =-\(\sum_{j,k=1}^n h^{jk} f_{x_j} f_t\)_{x_k}+\(\frac{1}{2}\sum_{j,k}^n h^{jk}f_{x_j}f_{x_k}\)_t+\frac{1}{2}\(\l_0 f^2+ w^2\)_t\\
\ns\ds= w_t(-L+\l_0)^{-1}\(qw+\sum_{j,k}^n q_1^k w_{x_k}+q_2w_t\).
\ea
\ee
Integrating it on $ \O $, using integration by parts, by (\ref{bca3}) and  the H${\rm \ddot{o}} $lder inequality, the Sobolev embedding theorem, we obtain
\begin{equation}\label{d13}
\ba{ll}
\ds \frac{dE(t)}{dt}\le 2C\int_\O\[ w_t (-L+\l_0)^{-1}\(qw+\sum_{k=1}^n q_1^k w_{x_k}+q_2 w_t+\l_0 w\) \]dx\\
\ns\ds\q\q~~  \le CrE(t),\q \forall~t\in [0,T],
\ea
\end{equation}
which can yield (\ref{abc20}).\q \endpf
\par Finally, we recall the following energy estimate.
\begin{lemma}\label{l2}
Let $ 0\le S_1<S_2<T_2<T_1\le T$ and $q(\cd)\in L^\infty (Q) $ and $ q_1^k(\cd)\in L^\i (Q) $  ($ k=1,\cd\cd\cd,n $), $ q_2(\cd)\in L^\i(Q) $. Then the weak  solution $ w(\cd)\in C([0,T];L^2(\O))\cap C^1([0,T];H^{-1}(\O)) $ of (\ref{aaa1}) satisfies
\begin{equation}\label{a26}
\int_{S_2}^{T_2}E(t)dt\le C(1+r^2 )\int_{S_1}^{T_1}|w(t,\cd)|^2_{L^2(\O)}dt.
\end{equation}
\end{lemma}
\par {\it Proof of Lemma \ref{l2}. }Denote\begin{equation}\label{d15}
\phi(t)=\left\{\ba{ll}
\ds 1,\q\q t\in [S_2,T_2],\\
\ns\ds 0, \q\q t\in \(-\infty, \frac{S_1+S_2}{2}\]\cup \[\frac{T_1+T_2}{2}, +\infty\).
\ea\right.
\end{equation}
Put $ g=(-L+\l_0)^{-1} w $ and by the following Pointwise identity,
\bel{bca1}
\ba{ll}
\ds \phi w (-L+\l_0)^{-1}\(w_{tt}-\sum_{j,k=1}^n (h^{jk}w_{x_j})_{x_k}\)\\
\ns\ds =(\phi w g_t)_t-\phi w_t g_t-\phi_t w g_t +\phi w^2-\l_o \phi w g\\
\ns\ds =(\phi w g_t)_t-\phi (-L+\l_0) g_t g_t -\phi_t (-L+\l_0) g g_t +\phi w^2-\l_0 \phi w g\\
\ns\ds =(\phi w g_t)_t +\sum_{k=1}^n \(\sum_{j=1}^n \phi h^{jk} g_{tx_j} g_t\)_{x_k}-\(\sum_{j,k=1}^n \phi h^{jk} g_{tx_j}g_{tx_k}\)-\l_0 \phi g_t^2\\
\ns\ds \q +\sum_{k=1}^n \(\sum_{j=1}^n \phi_t h^{jk} g_{x_j}g_{t}\)_{x_k}-\(\frac{1}{2}\sum_{j,k=1}^n \phi_t h^{jk} g_{x_j}g_{x_k}\)_t+\frac{1}{2} \phi_{tt}\(\sum_{j,k=1}^n h^{jk}g_{x_j}g_{x_k}\)\\
\ns\ds\q  -\(\frac{1}{2}\l_0 \phi_t g^2\)_t+\frac{1}{2} \l_0 \phi_{tt}g^2+\phi w^2-\l_0\phi w g,
\ea
\ee
integrating by part, we have
\begin{equation}\label{bca2}
\ba{ll}
\ns\ds \int_{S_1}^{T_1}\int_\O\phi w (-L+\l_0)^{-1}\(w_{tt}-\sum_{j,k=1}^n(h^{jk}w_{x_j})_{x_k}\) dxdt\\
\ns\ds =-\int_{S_1}^{T_1}\int_\O \sum_{j,k=1}^n\phi h^{j k}g_{t x_j}g_{t x_k} +\l_0 \phi g_t^2 dx dt\\
\ns\ds \q +\int_{S_1}^{T_1}\int_\O \frac{1}{2} \phi_{t t}\(\sum_{j,k=1}^n h^{j k}g_{x_j}g_{x_k}+\l_0 g^2 \)dx d t \\
\ns\ds \q -\int_{S_1}^{T_1}\int_\O\phi (w^2-\l_0 w g )dxdt.
\ea
\end{equation}
Thus, by Lemma \ref{lll}, Poincar$ {\rm\acute{e}} $ inequality and  equation (\ref{aaa1}), we get
\begin{equation}\label{d16}
\ba{ll}
\ds  C_1 h_0\int_{S_2}^{T_2}|w_t(t,\cd)|^2_{H^{-1}(\O)}dt\\
\ns\ds \le  \int_{S_1}^{T_1}\int_\O \sum_{j,k=1}^n\phi h^{jk}g_{tx_j}g_{tx_k} +\l_0 \phi g_t^2dxdt\\
\ns\ds  =-\int_{S_1}^{T_1}\int_\O \phi w (-L+\l_0)^{-1} (qw+\sum_{k=1}^n q_1^k w_{x_k}+q_2 w_t+\l_0 w) dxdt \\
\ns\ds \q +\frac{1}{2}\int_{S_1}^{T_1}\int_\O \phi_{tt}\(\sum_{j,k=1}^n h^{jk} g_{x_j}g_{x_k}+\l_0 g^2\)dxdt\\
\ns\ds \q+\int_{S_1}^{T_1}\int_\O \phi (w^2-\l_0 w g) dxdt\\
\ns\ds  \le C(1+r^2)\int_{S_1}^{T_1}|w(t,\cd)|^2_{L^2(\O)}dt+\frac{C_1 h_0}{2}\int_{S_1}^{T_1} \phi |w_t(t, \cd)|_{H^{-1}(\O)}^2dt.
\ea
\end{equation}
which implies the desired result immediately.\q\endpf

\section{Proof of Theorem \ref{theorem1}}

In this Section, we will give the proof of Theorem \ref{theorem1}.  We divide the proof into several steps.

\ss

{\bf Step 1.} Put
\begin{equation}\label{a27}
z(t,s,x)\=\int_s^tw(\tau, x)d\tau,\q \forall (t,s,x)\in (0,T)\t Q.
\end{equation}
Then, it is easy to see that $z$ satisfies the following ultra-hyperbolic equation:
\begin{equation}\label{a28}
\left\{\ba{ll}
\ds
z_{tt}+z_{ss}-\sum_{j,k=1}^n(h^{jk}z_{x_j})_{x_k}=F, &\ds\mbox{ in } (0,T)\t Q,\\
\ns\ds z=0,&\ds \mbox{ on }(0,T)\t \Sigma,
\ea\right.
\end{equation}
where
\bel{1207-f}
F\=\int_s^t q(\tau,x)z_t(\tau,s,x)+\sum_{k=1}^n q_1^k(\tau,x)z_{x_k,t}(\tau, s,x)+q_2(\tau,x)z_{tt}(\tau,s,x)d\tau.
\ee

For $\l>0$, we introduce the following weight functions:
\begin{equation}\label{a29}
\left\{\ba{ll}
\ns\ds \th(t,s,x)=e^{\ell(t,s,x)},\q \q\ell(t,s,x)=\l \phi(t,s,x),\\
\ns\ds \phi(t,s, x)=d(x)-\a (t-T/2)^2-\a (s-T/2)^2,
\ea\right.
\end{equation}
where $ \a\in (0,1) $ and  $ d(x) $ satisfies Condition \ref{cd1} and (\ref{a15}).

\ss

{\bf Step 2.}  Denote
\begin{equation}\label{a31}
\left\{\ba{ll}
\cQ\=(0,T)\t Q,\q ~~\q\q \cS\=(0,T)\t \G,\\
\ns\ds T_i\=\frac{T}{2}-\e_iT,\q\q\q\q T_i'\=\frac{T}{2}+\e_i T,\\
\ns\ds \cQ_i\=(T_i,T_i')\t (T_i,T_i')\t \O, \q\q \cQ_i'\=(T_i,T_i')\t(T_i,T_i')\t \o_1,\\
\ns\ds \cS_i\=(T_i,T_i')\t(T_i,T_i')\t \G,\q\q~ \cS_{0i}\=(T_i,T_i')\t(T_i,T_i')\t \G_0.
\ea
\right.
\end{equation}
where $ \o_1 $  is  some neighborhood of $ \G_0 $ satisfying $\o_0 \subset \o_1\subset \o$, and $ 0<\e_0<\e_1<\e_2<\e_3<\frac{1}{2} $.

\begin{figure}[!h]
\footnotesize
\begin{center}
\begin{tabular}{ccc}
\includegraphics[width=0.3\textwidth]{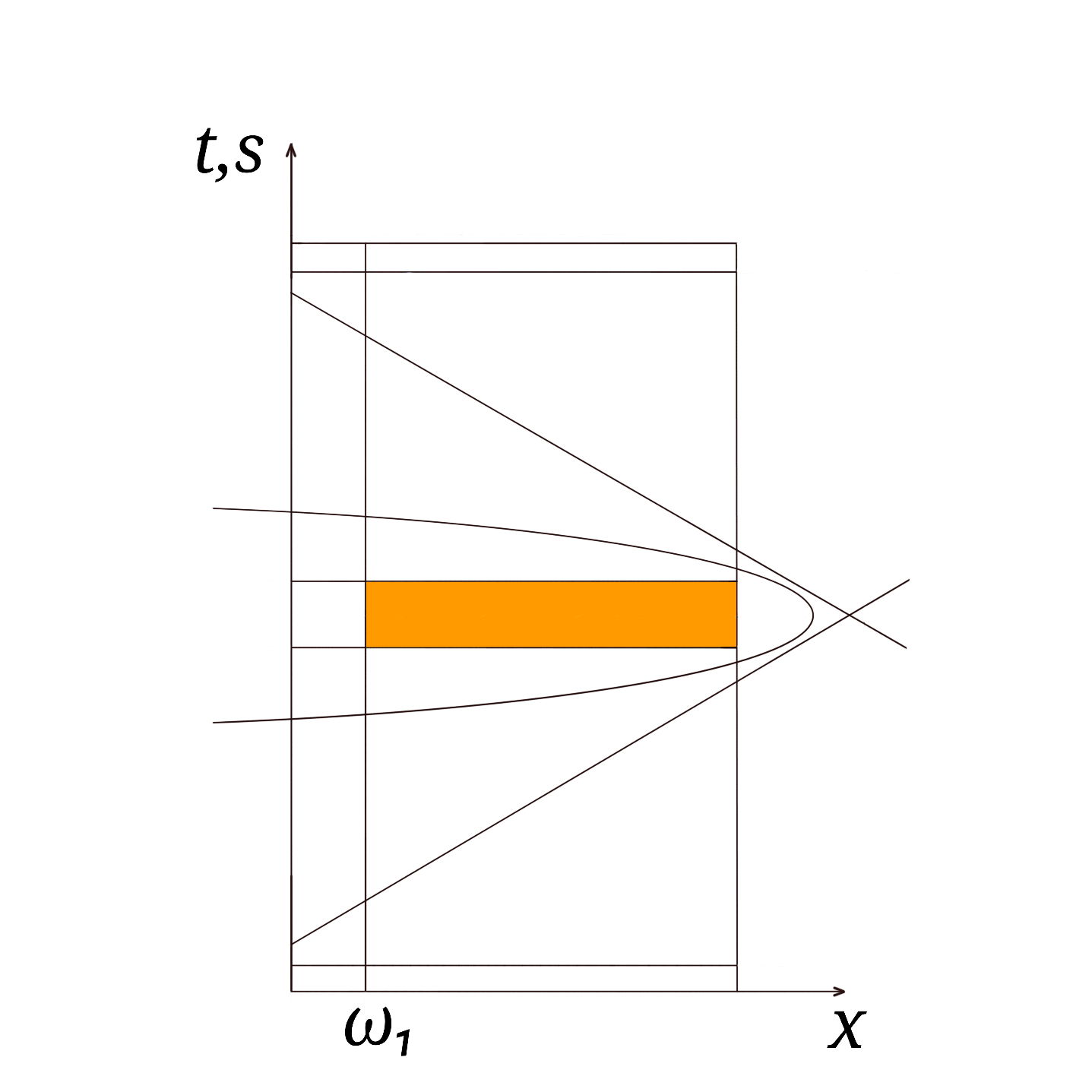}&
\hspace{-0.46cm}
\includegraphics[width=0.3\textwidth]{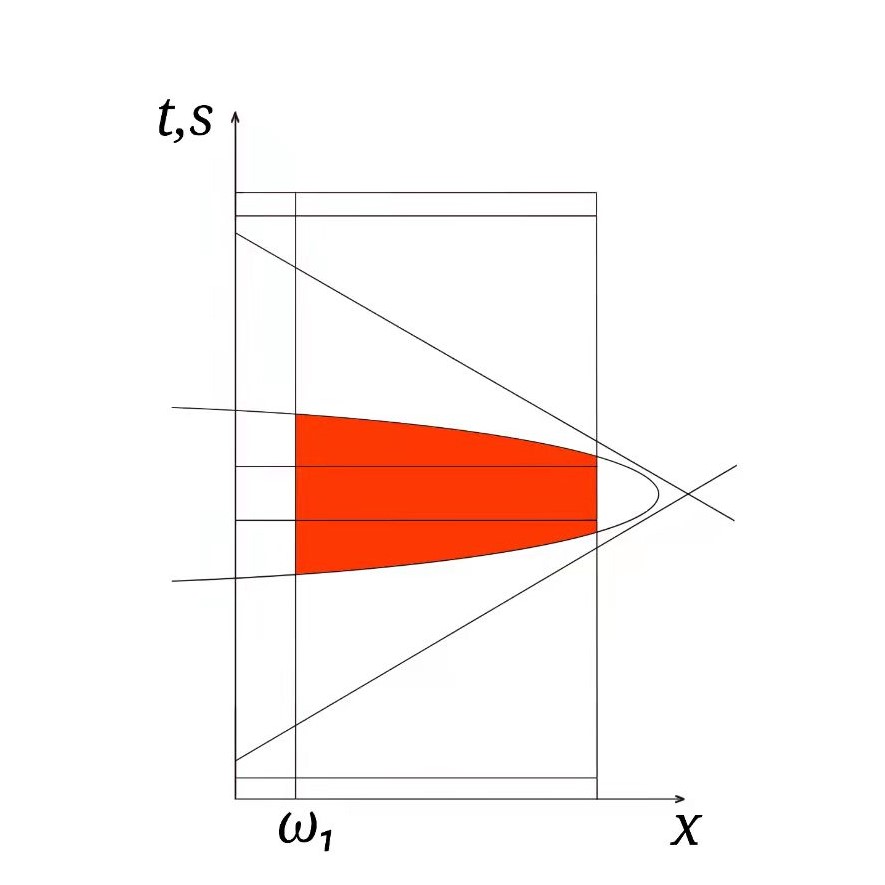}&
\hspace{-0.46cm}
\includegraphics[width=0.3\textwidth]{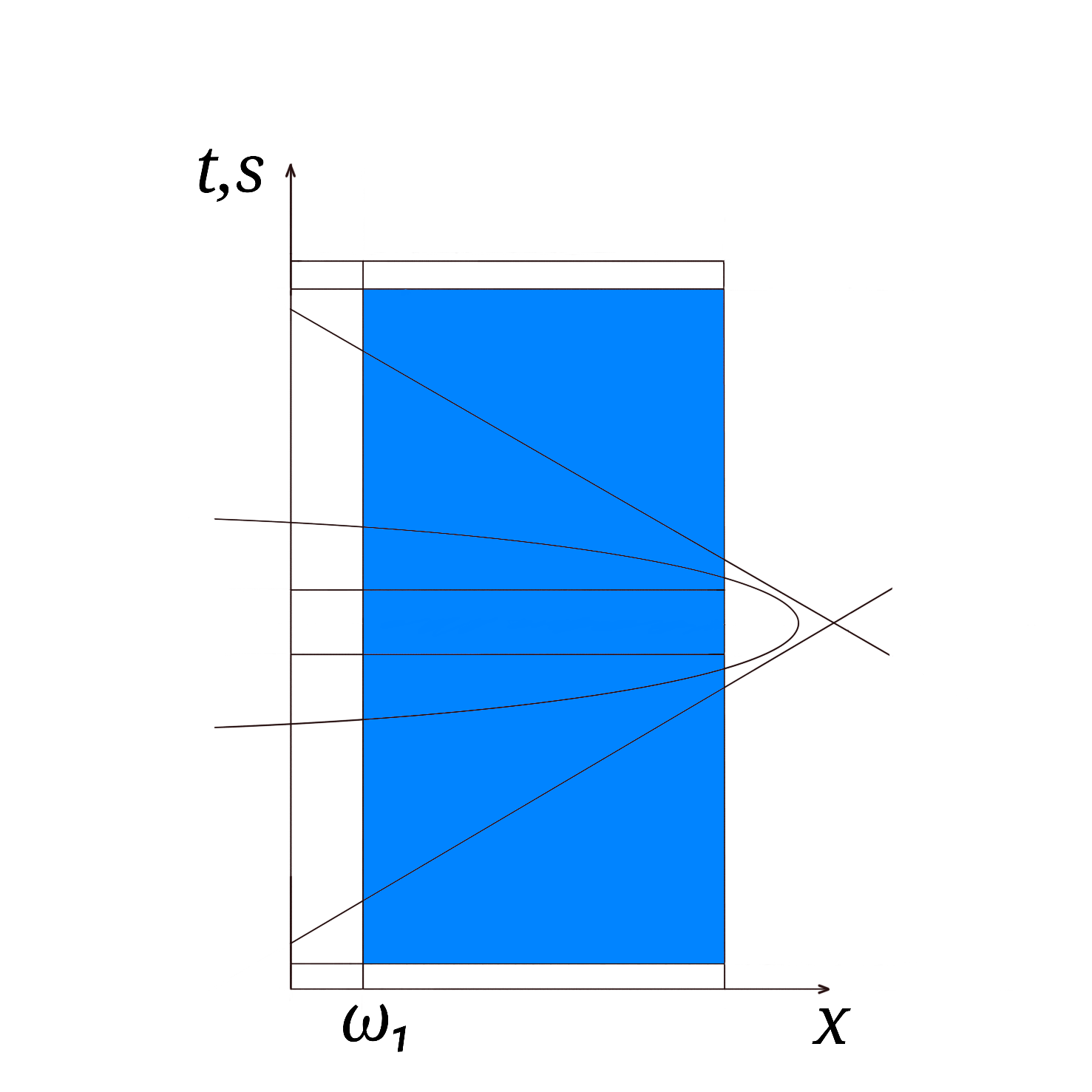}\\
(a)&(b)&(c)
\hspace{-0.46cm}
\end{tabular}
\end{center}
\caption{(a) represent $ \cQ_0\backslash \cQ_0'$, (b) represent $ \cQ(c)$, (c) represent $ \cQ_1\backslash \cQ_1'$.}\label{Fig3}
\end{figure}

Recalling (\ref{abc7}) the definitions of $R_0$ and $R_1$, since $ \ds \min_{x\in \oO} d(x)>0 $, we see that $ R_0>0 $.
We can choose a sufficiently small $ c\in (0,R_0) $ and an $ \a\in (0,1) $ close to 1 such that
\begin{equation}\label{aaaa1}
\min_{x\in \oO} d(x) >c^2,
\end{equation}
and
\begin{equation}\label{aaaaa1}
1-\frac{2c^2}{T^2}<\a<1.
\end{equation}
Hence, if  $ t=s=\frac{T}{2} $, for any $ x\in \cl{\O} $, we have
\begin{equation}\label{aa1}
\phi\(\frac{T}{2},\frac{T}{2},x\)=d(x)> c^2.
\end{equation}
For any $ b>0 $, define
\begin{equation}\label{a32}
\cQ(b)\=\Big\{(t,s,x)\in (-\infty,+\infty)\t (-\infty,+\infty)\t \{\O\backslash \o_1\}\big|~ \phi(t,s,x)>b^2\Big\}.
\end{equation}
Notice that for any  $x\in \cQ(c)  $, by (\ref{aaaaa1}), we have
\begin{equation}\label{aab2}
\ba{ll}
\ds d(x)-(t-T/2)^2-(s-T/2)^2\\
\ns\ds>\phi(t,s,x)+(\a-1)(t-T/2)^2+(\a-1)(s-T/2)^2\\
\ns\ds  >c^2+(\a-1)\frac{T^2}{2}>0.
\ea
\end{equation}
Noting that $T>2R_1$,  Choose $\o_1$ and $\o$ closely enough to make $ d(x)<T^2/4$ for any $x\in \O\backslash \o_1$. Now take $ \e_0  $ sufficiently small and $ \e_1\in (0,\frac{1}{2}) $ sufficiently close to $ \frac{1}{2} $ such that
\begin{equation}\label{a33}
\cQ_0\backslash \cQ_0'\subset \cQ(c)\subset\cD'\subset \cQ_1\backslash\cQ_1',
\end{equation}
where
\begin{equation}\label{a34}
\cD'\=\Big\{(t,s,x)\in \{\cQ\backslash (0,T)\t(0,T)\t \o_1\}\big|~ d(x)-(t-T/2)^2-(s-T/2)^2>0\Big\},
\end{equation}
and
\begin{equation}\label{aa2}
(T_0,T_0')\t (T_0,T_0')\t\O\subset \cD''\=\Big\{(t,s,x)\in \cQ \big|~ d(x)-(t-T/2)^2-(s-T/2)^2>0  \Big\}.
\end{equation}
We refer to  Figure \ref{Fig3} for the relations of  (\ref{a33}) in one dimensional case.
\\ Choosing $ \e $ small enough, we can still hold
\begin{equation}\label{a35}
\cQ_0\backslash \cQ_0'\subset\cQ(c+2\e)\subset\cQ(c+\e)\subset \cQ(c)\subset\cD' \subset \cQ_1\backslash\cQ_1'.
\end{equation}
Next, we introduce the following cut-off function:
\begin{equation}\label{a36}
\eta(t,s,x)=\left\{
\ba{ll}
\ds 1,& (t,s,x)\in \cQ(c+2\e)\backslash \{(0,T)\t (0,T)\t w_2\},\\
\ns\ds 0, & (t,s,x)\in \cQ_1\backslash \cQ(c+\e),
\ea\right.
\end{equation}
where $\o_2 $ is also some neighborhood of $\G_0$ satisfying $ \o_1\subset \o_2\subset \o $. Set $$ u(t,s,x)\=\eta (t,s,x)z(t,s,x). $$ It is easy to check that $u$ satisfies
\begin{equation}\label{a38}
\ba{ll}
\ds u_{tt}+u_{ss}-\sum_{j,k=1}^n(h^{jk}u_{x_j})_{x_k}\\
\ns\ds  ~\q=\eta F+\[\eta_{tt}+\eta_{ss}-\sum_{j,k=1}^n(h^{jk}\eta_{x_j})_{x_k}\]z+2\eta_{t}z_t+2\eta_{s}z_s -2\sum_{j,k=1}^nh^{jk}\eta_{x_j} z_{x_k},
\ea
\end{equation}
where $F$ is given by (\ref{1207-f}).

\ss

{\bf Step 3. } By (\ref{a29}), it is easy to see that
\bel{1207-l}\left\{\ba{ll}\ds
\ell_t= -2\l\a (t-T/2) , \q \ell_s=-2\l \a (s-T/2),\\
\ns\ds
\ell_{tt}=\ell_{ss}=-2\l \a ,\\
\ns\ds \ell_{x_j}=\l d_{x_j},~ \ell_{x_jx_k}=\l d_{x_jx_k}, \q j, k=1,...,n.
\ea\right.\ee
In the sequel, we denote $O(\l^k)$ by a function of order $ \l^k $ for large $ \l $. Recall Corollary \ref{Co4.1}, put
$$ \Psi=-\l \sum_{j,k=1}^n (h^{jk}d_{x_j})_{x_k}+2\l (1-\a). $$
By  (\ref{a8}),  (\ref{a15}) and (\ref{1207-l}), we have
\begin{equation}\label{baa1}
\ba{ll}
\ns\ds 2\[\ell_{tt}-\ell_{ss}+\sum_{j,k=1}^n (h^{jk}\ell_{x_j})_{x_k}+\Psi\]=4\l(1-\a),
\ea
\end{equation}
and
\bel{1207-cjk}\ba{ll}\ds
2\sum_{j,k=1}^n  c^{jk}v_{x_j}v_{x_k}\\
\ns\ds \ds=2\sum_{j,k=1}^n\Big\{\sum_{j',k'=1}^n\[2h^{jk'}(h^{j'k}\ell_{x_{j'}})_{x_{k'}}-(h^{jk}h^{j'k
	'}\ell_{x_{j'}})_{x_{k'}}\]\\
\ns\ds \qq\q\q\q+h^{jk}(\ell_{tt}+\ell_{ss}-\Psi )\Big\}v_{x_j}v_{x_k}\\
\ns\ds\ge 2\mu_0 \sum_{j,k=1}^n h^{jk}v_{x_j}v_{x_k}-8\l\a \sum_{j,k=1}^n h^{jk} v_{x_j}v_{x_k}-4\l(1-\a)\sum_{j,k=1}^n h^{jk}v_{x_j}v_{x_k}\\
\ns\ds\ge 4\l (1-\a)h_0|\n v|^2,

\ea\ee
Next, recall (\ref{c2}) for the definition of $A$ and $ B$, we have
\bel{1207-A}\ba{ll}\ds
A&\ds=\sum_{j,k=1}^n (h^{jk}\ell_{x_j}\ell_{x_k}-h_{x_j}^{jk}\ell_{x_k}-h^{jk}\ell_{x_jx_k})-\ell_t^2-\ell_s^2+\ell_{tt}+\ell_{ss}-\Psi\\
\ns&\ds=\l^2\[\sum_{j,k=1}^nh^{jk}d_{x_j}d_{x_k}-4a^2(t-T/2)^2-4\a^2(s-T/2)^2\]+O(\l).
\ea\ee
and
\begin{equation}\label{e3}
\ba{ll}
\ds
B =2A\[\Psi-\ell_{tt}-\ell_{ss}+\sum_{j,k=1}^n(h^{jk}\ell_{x_j})_{x_k}\]+2\(\sum_{j,k=1}^nh^{jk}\ell_{x_j}A_{x_k}-A_t\ell_t-A_s\ell_s\)\\
\ns\ds\q =4\l(1+\a)A+2\l^3\sum_{j,k=1}^nh^{jk}d_{x_j}\(\sum_{j',k'=1}^nh^{j'k'}d_{x_j'}d_{x_k'}\)_{x_k}\\
\ns\ds \q\q -32\l^3\a^3(t-\frac{T}{2})^2-32\l^3\a^3(s-\frac{T}{2})^2+O(\l^2)\\
\ns\ds \q \ge 4(3+\a)\l^3\[\sum_{j,k=1}^nh^{jk}d_{x_j}d_{x_k}-4\a^2(t-\frac{T}{2})^2-4\a^2(s-\frac{T}{2})^2\]+O(\l^2),
\ea
\end{equation}
where we used
\begin{equation}\label{ja1}
\ba{ll}
\ds \sum_{j,k=1}^nh^{jk}d_{x_j}\(\sum_{j',k'=1}^nh^{j'k'}d_{x_j'}d_{x_k'}\)_{x_k}\\
\ns\ds =\sum_{j,k=1}^n\sum_{j',k'=1}^n\[2h^{jk'}(h^{j'k}d_{x_j
	'})_{x_k'}-h_{x_k'}^{jk}h^{j'k'}d_{x_j'}\]d_{x_j}d_{x_k}\\
\ns\ds \ge \mu_0 \sum_{j,k=1}^n h^{jk}d_{x_j}d_{x_k} \ge 4\sum_{j,k=1}^n h^{jk}d_{x_j}d_{x_k}.
\ea
\end{equation}
Further, notice that
\bel{1207-mix}\ba{ll}\ds
\|-2\sum_{j,k=1}^nh^{jk}\Psi_{x_j}vv_{x_k}\|\le C\sum_{j,k=1}^n h^{jk} v_{x_j}v_{x_k}+O(\l^2) v^2.

\ea\ee
By (\ref{1207-cjk}), (\ref{e3}) and (\ref{1207-mix}), we conclude that there is a $\l_0>0$, for any $\l>\l_0$, we have
\begin{equation}\label{e2}
\ba{ll}
\ds \th^2\|u_{tt}+u_{ss}-\sum_{j,k=1}^n (h^{jk}u_{x_j})_{x_k}\|^2+2\div V+2 M_t+2N_s\\
\ns\ds \ge 2\l (1-\a)\(v_t^2+v_s^2+h_0|\n v|^2\)\\
\ns\ds\q+8(3+\a)\l^3\[d(x)-\a^2(t-\frac{T}{2})^2-\a^2(s-\frac{T}{2})^2\]v^2,
\ea
\end{equation}
Now,  integrate (\ref{e2}) on  $ \cQ_1 $,  we get
\begin{equation}\label{c7}
\ba{ll}
\ds \l\int_{\cQ(c+2\e)\backslash \{(0,T)\t (0,T)\t w_2\}} \th^2\[\(z_t^2+z_s^2+|\n z|^2\)+\l^2z^2\]dxdtds\\
\ns\ds\le  C\int_{\cQ_1}\th^2\|u_{tt}+u_{ss}-\sum_{j,k=1}^n (h^{jk}u_{x_j})_{x_k}\|^2dxdtds.\ea
\end{equation}

{\bf Step 4.}  Let us estimate ``$\ds \int_{\cQ_1}\th^2\|u_{tt}+u_{ss}-\sum_{j,k=1}^n (h^{jk}u_{x_j})_{x_k}\|^2dxdtds$".
Recall the definition of $ \eta $, noting that $ \eta\equiv 1 $ on $ \cQ(c+2\e)\backslash \{(0,T)\t (0,T)\t \o_2\} $, by (\ref{1207-f}) and (\ref{a38}), we have
\begin{equation}\label{a39}
\ba{ll}
\ds \int_{\cQ_1} \th^2\|u_{tt}+u_{ss}-\sum_{j,k=1}^n (h^{jk}u_{x_j})_{x_k}\|^2dxdtds\\
\ns\ds \le C\int_{\cQ_1}\th^2\|\int_s^tq(\tau,x)z_t(\tau,s,x)d\tau\|^2dxdtds\\
\ns\ds \q +C\int_{\cQ_1}\th^2\|\int_s^t \sum_{k=1}^n q_1^k(\tau,x)z_{x_k,t}(\tau, s,x)+q_2(\tau,x)z_{tt}(\tau,s,x)d \tau\|^2dxdtds \\
\ns\ds \q +Ce^{(c+2\e)^2\l}\int_{\cQ_1}\(z_t^2+z_s^2+|\n z|^2+z^2\)dxdtds\\

\ns\ds \q +C\int_{(T_1,T_1')\t (T_1,T_1')\t \o_2\cap \cQ(c+2\e)}\th^2\(z_t^2+z_s^2+|\n z|^2+z^2\)dxdtds.
\ea
\end{equation}
Notice that
\begin{equation}\label{b1}
\ba{ll}
\ds\int_{\cQ_1}\th^2\|\int_{T/2}^tz_t^2(\tau, s ,x)d \tau \|dtdsdx\\
\ns\ds  =\int_\O dx\int_{T_1}^{T_1'}ds\left\{\int_{T_1}^{T_1'}e^{\l\[d(x)-\a(t-T/2)^2-\a(s-T/2)^2\]}\|\int_{T/2}^tz_t^2(\tau,s,x)d\tau\|dt\right\}\\
\ns\ds  =\int_\O dx\int_{T_1}^{T_1'}ds\left\{\int_{T_1}^{T/2}dt\(e^{\l\[d(x)-\a(t-T/2)^2-\a(s-T/2)^2\]}\int_t^{T/2}z_t^2(\tau,s,x)d\tau\)\right.\\
\ns\ds \q+\left.\int_{T/2}^{T_1'}dt\(e^{\l\[ d(x)-\a(t-T/2)^2-\a(s-T/2)^2\]}\int_{T/2}^tz_t^2(\tau,s,x)d\tau\)\right\}\\
\ns\ds \le C\int_\O\int_{T_1}^{T_1'}ds\int_{T_1}^{T_1'}e^{\l\[d(x)-\a(\tau-T/2)^2-\a(s-T/2)^2\]}z_t^2(\tau,s,x)d\tau\\
\ns\ds \le C\int_{\cQ_1}\th^2z_t^2dxdtds.
\ea
\end{equation}
Similarly,
\begin{equation}\label{b2}
\int_{\cQ_1}\th^2\|\int_{T/2}^s z_s^2(t,\tau,x)d\tau\|dxdtds\le C\int_{\cQ_1}\th^2z_s^2dxdtds.
\end{equation}
Combining (\ref{b1}) and (\ref{b2}), noting that $z_t(\tau, s, x)=-z_s(t,\tau,x)$, we have
\bel{ha4}
\int_{\cQ_1}\th^2\|\int_s^t q(\tau,x)z_t(\tau,s,x)d\tau\|^2 dxdtds\le C r^2\int_{\cQ_1}\th^2\(z_t^2+z_s^2\)dxdtds.
\ee
Now, we consider the next term in (\ref{a39}). Using integration by parts, we have
\begin{equation}\label{bcc1}
\ba{ll}
\ds \int_{\cQ_1}\th^2\|\int_s^t\( \sum_{k=1}^n q_1^k(\tau,x)z_{tx_k}(\tau,s,x)+q_2(\tau,x)z_{tt}(\tau,s,x)\)d\tau\|^2dxdtds\\
\ns\ds \le \int_{\cQ_1}\th^2\|\sum_{k=1}^nq_1^k(t,x)z_{x_k}(t,s,x)\|^2+\th^2\|q_2(t,x)z_t(t,s,x)\|^2dxdtds\\
\ns\ds \q +\int_{\cQ_1}\th^2\|\int_s^t \sum_{k=1}^n q_{1,t}^k(\tau,x)z_{x_k}(\tau,s,x)+q_{2,t}(\tau,x)z_t(\tau,s,x)d\tau\|^2dxdtds.
\ea
\end{equation}
Similar to the argument in (\ref{b1}) and (\ref{b2}), we have
\begin{equation}\label{bcc2}
\ba{ll}
\ds \int_{\cQ_1}\th^2\|\int_s^t\( \sum_{k=1}^n q_1^k(\tau,x)z_{x_k,t}(\tau,s,x)+q_2(\tau,x)z_{tt}(\tau,s,x)\)d\tau\|^2dxdtds\\
\ns\ds \le Cr^2\int_{\cQ_1}\th^2 (|\n z|^2+z_t^2)dxdtds.
\ea
\end{equation}
Combining (\ref{a39}), (\ref{ha4}) and (\ref{bcc2}), we have
\begin{equation}\label{b3}
\ba{ll}
\ds \int_{\cQ_1}\th^2\| u_{tt}+u_{ss}-\sum_{j,k=1}^n( h^{jk}u_{x_j})_{x_k}\|^2dxdtds\\
\ns\ds \le C(1+r^2)e^{(c+2\e)^2\l}\int_{\cQ_1}\(z_t^2+z_s^2+|\n z|^2+z^2\)dxdtds\\
\ns\ds \q +C(1+r^2)\int_{(T_1,T_1')\t (T_1,T_1')\t \o_2\cap \cQ(c+2\e)}\th^2 \(z_t^2+z_s^2+|\n z|^2+z^2\)dxdtds\\
\ns\ds \q +Cr^2\int_{\cQ(c+2\e)\backslash \{(T_1,T_1')\t (T_1,T_1')\t \o_2\}}\th^2\(z_t^2+z_s^2+|\n z|^2\)dxdtds.
\ea
\end{equation}
Combining (\ref{c7}) and (\ref{b3}), when $ \frac{\l}{2}\ge Cr^2 $, we can yield
\begin{equation}\label{b4}
\ba{ll}
\ds \l\int_{\cQ(c+2\e)\backslash \{(0,T)\t (0,T)\t w_2\}} \th^2\[\(z_t^2+z_s^2+|\n z|^2\)+\l^2 z^2\]dxdtds\\
\ns\ds \le C(1+r^2)e^{\l(c+2\e)^2}\int_{\cQ_1}\(z_t^2+z_s^2+|\n z|^2+z^2\)dxdtds\\
\ns\ds \q +C(1+r^2)e^{C\l}\int_{(T_1,T_1')\t (T_1,T_1')\t \o_2\cap \cQ(c+2\e)}\(z_t^2+z_s^2+|\n z|^2+z^2\)dxdtds.
\ea
\end{equation}
Adding $ \ds \l\int_{(T_0,T_0')\t (T_0,T_0')\t \o_2}\th^2 (z_t^2+z_s^2)dxdtds $ on the each side of (\ref{b4}), we have
\begin{equation}\label{b5}
\ba{ll}
\ds\l e^{\l(c+2\e)^2}\int_{\cQ_0} \(z_t^2+z_s^2\)dxdtds\\
\ns\ds \le \l\int_{\cQ_0}\th^2\(z_t^2+z_s^2\)dxdtds\\
\ns\ds \le C(1+r^2)e^{\l (c+2\e)^2}\int_{\cQ_1}\(z_t^2+z_s^2+|\n z|^2+z^2\)dxdtds\\
\ns\ds \q +C(1+r^2)\l e^{C\l}\int_{(T_1,T_1')\t (T_1,T_1')\t \o_2\cap \cQ(c+2\e)}\(z_t^2+z_s^2+|\n z|^2+z^2\)dxdtds\\
\ns\ds \q +C \l\int_{(T_0,T_0')\t (T_0,T_0')\t \o_2}\th^2 \(z_t^2+z_s^2\)dxdtds.
\ea
\end{equation}

\ms

{\bf Step 5.} Now, we estimate $\ds  \int_{\cQ_1}|\n z|^2dxdtds $ and \\
$ \ds \int_{\{(T_1,T_1')\t (T_1,T_1')\t \o_2\}\cap \cQ(c+2\e)}|\n z|^2 dxdtds $ from the right-hand side of (\ref{b5}) by standard method which can be found, for example, in \cite{XZ1, XQX}. Set
\begin{equation}\label{b6}
\chi(t,s)=(t-T_2)(T_2'-t)(s-T_2)(T_2'-s),
\end{equation}
so that
\begin{equation}\label{b7}
\ba{ll}
\ds \int_{\cQ_2}\chi zFdxdtds\\
\ns\ds =\int_{\cQ_2}\chi z\( z_{tt}+z_{ss}-\sum_{j,k=1}^n (h^{jk}z_{x_j})_{x_k} \)dxdtds\\
\ns\ds \ge -\int_{\cQ_2}\[z_t(\chi_tz+\chi z_t)+z_s(\chi_sz+\chi z_s)\]dxdtds+h_0\int_{\cQ_2}\chi|\n z|^2dxdtds.
\ea
\ee
and
\bel{ha5}
\ba{ll}
\ds \int_{\cQ_2}\chi z Fdxdtds\\
\ns\ds=\int_{\cQ_2} \chi z\(\int_s^t q(\tau,x)z_t(\tau,s,x)d\tau\)dxdtds\\
\ns\ds  \q +\int_{\cQ_2}\chi z\[\(\int_s^t -\pa_tq_2(\tau,x)z_t(\tau,s,x)d\tau\)+ q_2(t,x)z_t(t,s,x)\]dxdtds\\
\ns\ds \q +\int_{\cQ_2} \chi z\( \int_s^t \sum_{k=1}^n q_1^k(t,x)z_{tx_k}(\tau,s,x)d\tau\)dxdtds,
\ea
\ee
Combining with
\bel{haha1}
\ba{ll}
\ds \int_{\cQ_2}\chi z\(\int_s^t \sum_{k=1}^n  q_1^k(\tau, x) z_{tx_k}(\tau,s,x)d\tau \)dxdtds\\
\ns\ds =\int_{\cQ_2}\int_s^t\sum_{k=1}^n\(\chi(t,s)z(t,s,x) q_1^k(\tau,x)z_t(\tau,s,x)\)_{x_k}d\tau dxdtds\\
\ns\ds \q -\int_{\cQ_2}\int_s^t\[\chi(t,s)z_t(\tau,s,x)\(\sum_{k=1}^n z_{x_k}(t,s,x)q_1^k(\tau,x)\)\]d\tau dxdtds\\
\ns\ds \q -\int_{\cQ_2}\int_s^t\[\chi(t,s)z(t,s,x)z_t(\tau,s,x)\(\sum_{k=1}^n\pa_{x_k} q_1^k(\tau,x)\)\]d\tau dxdtds\\
\\
\ns\ds \le Cr^2\int_{\cQ_2}\(z^2+z_t^2 \)dx dt ds+c_0\int_{\cQ_2}\chi(t,s)|\n z|^2dxdtds,
\ea
\ee
by choosing $c_0$ sufficiently small, we get
\begin{equation}\label{b8}
\int_{\cQ_1}|\n z|^2dxdtds\le C(1+r^2)\int_{\cQ_2}\(z_t^2+z_s^2+z^2\)dxdtds.
\end{equation}
Similarly, recall $ \cD''=\Big\{(t,s,x)\in \cQ\big|~ d(x)-(t-T/2)^2-(s-T/2)^2>0\Big\} $ and notice that
\begin{equation}\label{ff1}
\big\{(T_1,T_1')\t (T_1.T_1')\t\o_2\big\} \cap  \cQ(c+2\e)\subset \cD''\cap \big\{(0,T)\t (0,T)\t \{\o\backslash \o_0\}\big\}.
\end{equation}
Putting
\begin{equation}\label{b9}
\varrho(t,s,x)=\left\{\ba{ll}
\ds 1,& (t,s,x)\in (T_1,T_1')\t (T_1.T_1')\t\o_2 \cap  \cQ(c+2\e),\\
\ns\ds 0,&(t,s,x)\in \cQ\backslash\big\{\cD''\cap (0,T)\t (0,T)\t \{\o\backslash \o_0\}\big\},
\ea\right.
\end{equation}
we have
\begin{equation}\label{b10}
\ba{ll}
\ds \int_{\cQ}\varrho zF dxdtds\\
\ns\ds =\int_{\cQ}\varrho z\( z_{tt}+z_{ss}-\sum_{j,k=1}^n (h^{jk}z_{x_j})_{x_k} \)dxdtds\\
\ns\ds \ge -\int_{\cQ}\[z(\varrho_tz+\varrho z_t)+z_s(\varrho_sz+\varrho z_s)\]dxdtds+h_0\int_{\cQ}\varrho |\n z|^2dxdtds\\
\ns\ds \q +\int_\cQ\frac{1}{2} z^2\sum_{j,k=1}^n(h^{jk}\varrho_{x_j})_{x_k} dxdtds,
\ea
\end{equation}
and
\bel{hahaha1}
\ba{ll}
\ds \int_{\cQ}\varrho z Fdxdtds\\
\ns\ds=\int_{\cQ} \varrho  z\(\int_s^t q(\tau,x)z_t(\tau,s,x)d\tau\)dxdtds\\
\ns\ds  \q +\int_{\cQ}\varrho  z\[\(\int_s^t \pa_tq_2(\tau,x)z_t(\tau,s,x)d\tau\)+q_2(t,x)z_t(t,s,x)\]dxdtds\\
\ns\ds \q +\int_{\cQ} \varrho z\( \int_s^t \sum_{k=1}^n  q_1^k(\tau,x)z_{tx_k}(\tau,s,x)d\tau \)dxdtds.
\ea
\ee
Combining (\ref{b9}) and (\ref{hahaha1}) with
\bel{haha2}
\ba{ll}
\ds \int_{\cQ}\varrho z\(\int_s^t \sum_{k=1}^n  q_1^k(\tau, x) z_{tx_k}(\tau,s,x)d\tau \)dxdtds\\
\ns\ds =\int_{\cQ}\int_s^t\sum_{k=1}^n\(\varrho(t,s,x)z(t,s,x) q_1^k(\tau,x)z_t(\tau,s,x)\)_{x_k}d\tau dxdtds\\
\ns\ds \q -\int_{\cQ}\int_s^t\[\(\sum_{k=1}^n z_{x_k}(t,s,x)q_1^k(\tau,x)\)\varrho(t,s,x)z_t(\tau,s,x)\]d\tau dxdtds\\
\ns\ds \q -\int_{\cQ}\int_s^t\[\varrho (t,s,x)z(t,s,x)z_t(\tau,s,x)\(\sum_{k=1}^n\pa_{x_k} q_1^k(\tau,x)\)\]d\tau dxdtds\\
\ns\ds \q -\int_{\cQ}\int_s^t\[\(\sum_{k=1}^n \varrho_{x_k}(t,s,x) q_1^k(\tau,x)\)z(t,s,x)z_t(\tau,s,x)\]d \tau dxdtds\\
\ns\ds \le Cr^2\int_{\cQ }\(z^2+z_t^2\)dx dt ds+c_1\int_{\cQ}\varrho|\n z|^2dxdtds,
\ea
\ee
by choosing $c_1$ sufficiently small, we yield
\begin{equation}\label{b11}
\ba{ll}
\ds \int_{\big\{(T_1,T_1')\t (T_1,T_1')\t \o_2\big\}\cap \cQ(c+2\e)}|\n z|^2dxdtds\\
\ns\ds \le C(1+r^2)\int_{\cD''\cap\big\{ (0,T)\t (0,T)\t \{\o\backslash \o_0\}\big\}}z_s^2+z_t^2+z^2dxdtds.
\ea
\end{equation}
Combining (\ref{b5}), (\ref{b8}) and (\ref{b11}), we get
\begin{equation}\label{b12}
\ba{ll}
\ds \l\int_{\cQ_0}(z_t^2+z_s^2)dxdtds\\
\ns\ds \le C(1+r^4)\int_{\cQ_2}(z_t^2+z_s^2+z^2)dxdtds\\
\ns\ds \q +C(1+r^4)\l e^{C\l}\int_{\cD''\cap\big\{ (0,T)\t (0,T)\t\{ \o\backslash \o_0\}\big\}}(z_s^2+z_t^2+z^2)dxdtds\\
\ns\ds \q +\l e^{C\l}\int_{(T_0,T_0')\t (T_0,T_0')\t \o_2}(z_t^2+z_s^2)dxdtds
\ea
\end{equation}

{\bf Step 6.} Let us return to the function $ w $. Recall that $ \ds z(t,s,x)=\int_s^t w(\tau,x)dxd\tau $. By (\ref{b12}), we have
\begin{equation}\label{b13}
\ba{ll}
\ds \l \int_{T_0}^{T_0'}\int_\O w^2dxdt\\
\ns\ds \le C(1+r^4)\l \int_{\cQ_2}\[w^2(t,x)+\(\int_s^t w(\tau)d\tau\)^2\]dxdtds\\
\ns\ds \q +C\l e^{C\l}(1+r^4)\int_{\cD''\cap \big\{(0,T)\t (0,T)\t\{ \o\backslash \o_0\}\big\}}\(w^2(t,x)+\(\int_s^tw^2(\tau,x)d\tau\)\)dxdtds\\
\ns\ds \q +\l e^{C\l}\int_{(T_0,T_0')\t \o_2}w^2(t,x)dxdt\\
\ns\ds  \le C(1+r^4)\int_Q w^2(t,x)dxdt+C(1+r^4)\l e^{C\l}\int_{\cD\cap\big\{(0,T)\t \{\o\backslash \o_0\}\big\}} w^2(t,x)dxdt\\
\ns\ds \q +\l e^{C\l}\int_{(T_0,T_0')\t \o_2}w^2(t,x)dxdt,
\ea
\end{equation}
where we used
\begin{equation}\label{g2}
\ba{ll}
\ds \int_{\cD''\cap \big\{(0,T)\t (0,T)\t \{\o\backslash \o_0\}\big\}}\(\int_s^tw(\tau,x)d\tau \)^2dxdtds\\
\ns\ds\le \int_{\cD''\cap \big\{(0,T)\t (0,T)\t \{\o\backslash \o_0\}\big\}}\int_{\frac{T}{2}-\sqrt{d(x)}}^{\frac{T}{2}+\sqrt{d(x)}}w^2(\tau,x)d\tau dxdtds\\
\ns\ds  \le C\int_{\cD\cap\big\{(0,T)\t \{\o\backslash\o_0\}\big\}}w^2(t,x)dxdt,
\ea
\end{equation}
where $ \cD\=\Big\{(t,x)\in Q\big|~ d(x)-(t-T/2)^2>0\Big\} $ defined in (\ref{abc8}). Recall that $ \o=\cO_{\d}(\G_0) $ and $ \o_0=\cO_{\d_0}(\G_0) $, let $\e_0<\d_1$ so that
\begin{equation}\label{aa3}
\cD\cap \Big\{(0,T)\t \big\{\o\backslash\o_0\big\}\Big\}\cup \Big\{(T_0,T_0')\t \o_2\Big\}\subset \cK.
\end{equation}
\par On the other hand, set $ S_0\in (T_0,\frac{1}{2}T) $ and $ S_0'\in (\frac{1}{2}T,T_0') $. Applying lemma \ref{lemma31} and lemma \ref{l2}, we immediately get
\begin{equation}\label{b14}
\l \int_{S_0}^{S_0'}E(t)dt\le C\l e^{C\l}(1+r^6)\int_{\cK }w^2(t,x)dxdt+C(1+r^6)\int_0^TE(t)dt.
\end{equation}
Combining with energy estimate (\ref{abc20}), we get
\begin{equation}\label{aabb2}
\(C\l e^{Cr}-C(1+r^6)e^{Cr}\)E(0)\le C\l e^{C\l}(1+r^6)\int_{\cK} w^2(t,x)dxdt.
\end{equation}
So that, we have
\begin{equation}\label{b15}
E(0)\le \l e^{C\l}(1+r^6)\int_{\cK } w^2(t,x)dxdt,
\end{equation}
when $ \l\ge C(1+r^6)e^{Cr} $,
which yields (\ref{abc10}).\q\endpf
\section{Proof of Theorem \ref{theorem2}}

In this section, we will give the proof of Theorem \ref{theorem2}.

{\bf Proof. } Denote
\begin{equation}\label{f5}
\cQ_\zeta(b)\=\Big\{(t,s,x)\in (-\infty,+\infty)\t(-\infty,+\infty)\t \{\O\backslash\o_1\}\big|~\phi(t,s,x+\zeta)>b^2\Big\}.
\end{equation}
Notice that when $ t=\frac{T}{2} $, we have
\begin{equation}\label{aab1}
d(x)+d(x+\zeta)>0,\q \forall x\in \O.
\end{equation}
Set $\ds R_0\= \min_{x\in\oO} \Big\{\frac{1}{2} \(d(x)+d(x+\zeta)\)\Big\} $. Recall the note in (\ref{a31}) and (\ref{a32}), similar to the argument in (\ref{aaaa1}), (\ref{aaaaa1}) and (\ref{aa1}), we can still choose a sufficiently small $ c\in(0,R_0) $ and $ \a\in (0,1) $ close to $ 1 $ and take $ \e_0 $ small enough,  $ \e_1\in (0,\frac{1}{2}) $ sufficiently close to $ \frac{1}{2} $ such that
\begin{equation}\label{f6}
\cQ_0\backslash \cQ_0'\subset \cQ(c)\cup \cQ_\zeta(c)\subset \cD'\cup\tilde{\cD'}\subset \cQ_1\backslash \cQ_1',
\end{equation}
where $ \cD' $ is defined in (\ref{a34}) and
\begin{equation}\label{f9}
\ba{ll}
\tilde{\cD'}\=\Big\{(t,s,x)\in \cQ\backslash \{(0,T)\t (0,T)\t \o_1\}\big|~ d(x+\zeta )-(t-T/2)^2-(s-T/2)^2>0\Big\}.
\ea
\end{equation}
Choosing $ \e $ small enough, we have the following relations similar to (\ref{a35}),
\begin{equation}\label{f10}
\cQ_0\backslash \cQ_0'\subset \cQ(c+2\e)\cup \cQ_\zeta(c+2\e)\subset \cQ(c+\e)\cup \cQ_\zeta(c+\e)\subset \cQ(c)\cup \cQ_\zeta(c),
\end{equation}
and put $ \o_1\subset\o_2\subset \o $, by the same way we yield (\ref{b4}), we have
\begin{equation}\label{f11}
\ba{ll}
\ds \int_{\cQ_\zeta(c+2\e)\backslash \{(0,T)\t (0,T)\t w_2\}} \th^2(t,s,x+\zeta )\(\l(z_t^2+z_s^2+|\n z|^2)+\l^3 z^2\)dxdtds\\
\ns\ds \le C(1+r^2)e^{\l(c+2\e)^2}\int_{\cQ_1}\(z_t^2+z_s^2+|\n z|^2+z^2\)dxdtds\\
\ns\ds \q +C(1+r^2)e^{C\l}\int_{\{(T_1,T_1')\t (T_1,T_1')\t \o_2\}\cap \cQ_\zeta (c+2\e)}\(z_t^2+z_s^2+|\n z|^2+z^2\)dxdtds.
\ea
\end{equation}
By the definition of $ \th $ and $ Q_\zeta(c+2\e) $, we have
\begin{equation}\label{f12}
\ba{ll}
\ds e^{\l (c+2\e)^2}\int_{\cQ_\zeta (c+2\e)\backslash \{(0,T)\t (0,T)\t w_2\}} \(\l(z_t^2+z_s^2+|\n z|^2)+\l^3 z^2\)dxdtds\\
\ns\ds \le C(1+r^2)e^{\l(c+2\e)^2}\int_{\cQ_1}\(z_t^2+z_s^2+|\n z|^2+z^2\)dxdtds\\
\ns\ds \q +C(1+r^2)e^{C\l}\int_{\{(T_1,T_1')\t (T_1,T_1')\t \o_2\}\cap \cQ_\zeta (c+2\e)}\(z_t^2+z_s^2+|\n z|^2+z^2\)dxdtds.
\ea
\end{equation}
Combining with (\ref{b4}) and add $ \ds \l \int_{(T_0,T_0')\t (T_0,T_0')\t \o_2}\(z_t^2+z_s^2\)dxdtds $, we have
\begin{equation}\label{f13}
\ba{ll}
\ds \int_{\{\cQ(c+2\e)\cup \cQ_\zeta(c+2\e)\}\backslash\{(0,T)\t(0,T)\t \o_2\}}\(\l(z_t^2+z_s^2+|\n z|^2)+\l^3 z^2\)dxdtds\\
\ns\ds \q +\l \int_{(T_0,T_0')\t (T_0,T_0')\t \o_2}\(z_t^2+z_s^2\)dxdtds\\
\ns\ds \le C(1+r^2)\int_{\cQ_1}\(z_t^2+z_s^2+|\n z|^2+z^2\)dxdt\\
\ns\ds \q +C(1+r^2)\l e^{C\l}\int_{\{(T_1,T_1')\t (T_1,T_1')\t \o_2\}\cap \{\cQ(c+2\e)\cup\cQ_\zeta(c+2\e)\}}\(z_t^2+z_s^2+|\n z|^2+z^2\)dxdtds\\
\ns\ds \q +\l \int_{(T_0,T_0')\t (T_0,T_0')\t \o_2}\(z_t^2+z_s^2\)dxdtds.
\ea
\end{equation}
Recall that $ \cQ_0\backslash\cQ_0'\subset \cQ(c+2\e)\cup\cQ_\zeta(c+2\e) $, so that
\begin{equation}\label{f14}
\ba{ll}
\ds \int_{\cQ_0}\l(z_t^2+z_s^2)dxdtds\\
\ns\ds \le C(1+r^2)\int_{\cQ_1}(z_t^2+z_s^2+|\n z|^2+z^2)dxdt\\
\ns\ds \q +C(1+r^2)\l e^{C\l}\int_{\{(T_1,T_1')\t (T_1,T_1')\t \o_2\}\cap \{\cQ(c+2\e)\cup\cQ_\zeta(c+2\e)\}}(z_t^2+z_s^2+|\n z|^2+z^2)dxdtds\\
\ns\ds \q +\l \int_{(T_0,T_0')\t (T_0,T_0')\t \o_2}(z_t^2+z_s^2)dxdtds.
\ea
\end{equation}
Notice that
\begin{equation}\label{aa8}
\ba{ll}
\ns\ds (T_1,T_1')\t (T_1,T_1')\t \o_2\cap \{\cQ(c+2\e)\cup\cQ_\zeta(c+2\e)\}\\
\ns\ds \q\q\qq\subset \big\{\cD'\cup\tilde{\cD}'\big\}\cap \big\{(0,T)\t (0,T)\t\o\backslash\o_0\big\}. \ea
\end{equation}
By the same process of {\bf Step 5} in proof of Theorem \ref{theorem1}, we get
\begin{equation}\label{aa9}
\ba{ll}
\ds \int_{\cQ_0}(\l(z_t^2+z_s^2)dxdtds\\
\ns\ds \q\le C(1+r^4)\int_{\cQ_2}(z_t^2+z_s^2+z^2)dxdt\\
\ns\ds \qq +C(1+r^4)\l e^{C\l}\int_{\big\{\cD'\cup\tilde{\cD}'\big\}\cap \big\{(0,T)\t (0,T)\t\{\o\backslash\o_0\}\big\}}(z_t^2+z_s^2+z^2)dxdsdt\\
\ns\ds \qq +\l \int_{(T_0,T_0')\t (T_0,T_0')\t \o_2}(z_t^2+z_s^2)dxdtds.
\ea
\end{equation}
Notice that
\begin{equation}\label{aa10}
(T_0,T_0')\t \o_2\cup\Big\{\big\{\cD\cup\cD_\zeta \big\}\cap \big\{ (0,T)\t\o\backslash\o_0\big\}\Big\}\subset \cK\cup\cK_\zeta\subset W.
\end{equation}
using the same argument from {\bf Step 6} of the proof of Theorem
\ref{theorem1}, we can yield (\ref{abc14}).
\endpf
\section{Appendix A }

{\bf  Proof of Corollary \ref{Co4.1}. } Recall that $ \th=e^\ell $ and $ v=\th z $. Some elementary calculations yield that
\begin{equation}\label{bc1}
\ba{ll}
\ds \th\(z_{tt}+z_{ss}-\sum_{j,k=1}^n(h^{jk}z_{x_j})_{x_k}\)\\
\ns \ds =v_{tt}+v_{ss}-\sum_{j,k=1}^n(h^{jk}v_{x_j})_{x_k}-2\ell_tv_t-2\ell_s v_s+2\sum_{j,k=1}^n h^{jk}\ell_{x_j} v_{x_k}\\
\ns\ds \q+\ell_t^2v+\ell_s^2 v -\sum_{j,k=1}^nh^{jk}\ell_{x_j}\ell_{x_k}v-\ell_{tt}v-\ell_{ss}v+\sum_{j,k=1}^n(h^{jk}\ell_{x_j})_{x_k} v\\
\ns\ds   =I_1+I_2,
\ea
\end{equation}
where
\begin{equation}\label{bc2}
\left\{
\ba{ll}
\ns\ds I_1\= v_{tt}+v_{ss}-\sum_{j,k=1}^n(h^{jk}v_{x_j})_{x_k}-Av,\\
\ns\ds I_2\= -2\ell_t v_t-2\ell_s v_s+2\sum_{j,k=1}^n h^{jk}\ell_{x_j} v_{x_k}-\Psi v,\\
\ns\ds A\=-\ell_t^2-\ell_s^2+\sum_{j,k=1}^nh^{jk}\ell_{x_j}\ell_{x_k}+\ell_{tt}+\ell_{ss}-\sum_{j,k=1}^n (h^{jk}\ell_{x_j})_{x_k}-\Psi.
\ea\right.
\end{equation}
Then
\begin{equation}\label{bc3}
\th^2\|z_{tt}+z_{ss}-\sum_{j,k=1}^n(h^{jk}z_{x_j})_{x_k}\|^2= |I_1|^2+|I_2|^2+2I_1I_2\ge 2I_1I_2.
\end{equation}
\par Let us compute $ 2I_1I_2 $. Denote the terms in the right-hand sides of $ I_1 $ and $ I_2 $ by $ I_1^j~(j=1,2,3,4) $ and $ I_2^k~(k=1,2,3,4) $, respectively. Then
\begin{equation}\label{bc4}
\ba{ll}
\ds I_1^1I_2&\ds=\(-\ell_tv_t^2-2\ell_sv_sv_t+2\sum_{j,k=1}^nh^{jk}\ell_{x_j}v_{x_k}v_t-\Psi vv_t \)_t\\
\ns&\ds \q+\ell_{tt}v_t^2+2\ell_{st}v_sv_t+(\ell_s v_t^2)_s-\ell_{ss}v_t^2-2\sum_{j,k=1}^nh^{jk}\ell_{tx_j}v_{x_k}v_t\\
\ns&\ds \q-\div\(\sum_{j=1}^nh^{jk}\ell_{x_j}v_t^2\)+\sum_{j,k=1}^n(h^{jk}\ell_{x_j})_{x_k}v_t^2+\Psi v_t^2,
\ea
\end{equation}
and
\begin{equation}\label{bc5}
\ba{ll}
\ds I_1^2I_2&\ds=\(-\ell_sv_s^2-2\ell_tv_sv_t+2\sum_{j,k=1}^nh^{jk}\ell_{x_j}v_{x_k}v_s-\Psi vv_s \)_s\\
\ns&\ds\q+\ell_{s}v_s^2+2\ell_{st}v_sv_t+(\ell_t v_s^2)_t-\ell_{tt}v_s^2-2\sum_{j,k=1}^nh^{jk}\ell_{sx_j}v_{x_k}v_s\\
\ns&\ds\q-\div\(\sum_{j=1}^nh^{jk}\ell_{x_j}v_s^2\)+\sum_{j,k=1}^n(h^{jk}\ell_{x_j})_{x_k}v_s^2+\Psi v_s^2.
\ea
\end{equation}
For $ I_1^3I_2 $, we have
\begin{equation}\label{bc6}
\ba{ll}
\ds I_1^3I_2=-\sum_{j,k=1}^n(h^{jk}v_{x_j})_{x_k}\(-2\ell_tv_t-2\ell_sv_s+2\sum_{j,k=1}^nh^{jk}\ell_{x_j}v_{x_k}-\Psi v\)\\
\ns\ds ~~~~\,\,\,\,=\div\(\sum_{j,k=1}^n2h^{jk}v_{x_j}\ell_t v\)-\sum_{j,k=1}^n2h^{jk}v_{x_j}\ell_{tx_k}v_t-\(\sum_{j,k=1}^nh^{jk}\ell_t v_{x_j}v_{x_k}\)_t\\
\ns\ds \q\qq+\sum_{j,k=1}^n h^{jk}\ell_{tt}v_{x_j}v_{x_k}+\div\(\sum_{j,k=1}^n2h^{jk}v_{x_j}\ell_s v_s\)-\sum_{j,k=1}^n2h^{jk}v_{x_j}\ell_{sx_k}v_s\\
\ns\ds\q\qq-\(\sum_{j,k=1}^nh^{jk}\ell_s v_{x_j}v_{x_k}\)_s+\sum_{j,k=1}^n h^{jk}\ell_{ss}v_{x_j}v_{x_k}\\
\ns\ds \q\qq-\div\(2\sum_{j,j',k'=1}^nh^{jk}h^{j'k'}\ell_{x_{j'}}v_{x_j}v_{x_{k'}}\)+2\sum_{j,k,j',k'=1}^nh^{jk'}(h^{j'k}\ell_{x_{j'}})_{x_{k'}}v_{x_j}v_{x_k}\\
\ns\ds \q\qq+\div\(\sum_{j,j',k'=1}^nh^{jk}h^{j'k'}\ell_{x_j}v_{x_{j'}}v_{x_{k'}}\)-\sum_{j,k,j',k'=1}^n(h^{jk}h^{j'k'}\ell_{x_{j'}})_{x_{k'}}v_{x_j}v_{x_k}\\
\ns\ds \q\qq+\div\(\sum_{j,k=1}^nh^{jk}v_{x_j}\Psi v\)-\sum_{j,k=1}^nh^{jk}\Psi v_{x_j}v_{x_k}-\sum_{j,k=1}^nh^{jk}\Psi_{x_j}vv_{x_k}.
\ea
\end{equation}
For $I_1^4I_2$, we have
\begin{equation}\label{bc7}
\ba{ll}
\ds I_1^4I_2=-Av\(-2\ell_tv_t-2\ell_sv_s+2\sum_{j,k=1}^nh^{jk}\ell_{x_j}v_{x_k}-\Psi v\)\\
\ns\ds\qq =(A\ell_tv^2)_t-(A\ell_t)_tv^2+(A\ell_sv^2)_s-(A\ell_s)_sv^2\\
\ns\ds \q\qq-\div\(\sum_{j,k=1} h^{jk}A\ell_{x_j}v^2\)+\sum_{j,k=1}^n(Ah^{jk}\ell_{x_j})_{x_k}v^2+A\Psi v^2.
\ea
\end{equation}
Finally, combining  (\ref{bc3})--(\ref{bc7}), we complete the proof of Corollary \ref{Co4.1}. \endpf

\end{document}